\newcommand{\E}[0]{{\sf E}}

\documentclass[letterpaper,11pt]{article}

\usepackage{amsthm}
\usepackage{amsmath}

\usepackage{amsfonts}
\usepackage{amscd,amsthm,amsfonts,amsopn,amssymb,verbatim}

\newtheorem{thm}{Theorem}[section]

\newtheorem{lemma}[thm]{Lemma}
\newtheorem{prop}[thm]{Proposition}
\newtheorem{cor}[thm]{Corollary}

\newtheorem{conj}[thm]{Conjecture}

\newtheorem{question}[thm]{Question}

\newcommand{\beq}[1]{\begin{equation}\label{#1}}
\newcommand{\enq}[0]{\end{equation}}

\newcommand{\bn}[0]{\bigskip\noindent}
\newcommand{\mn}[0]{\medskip\noindent}
\newcommand{\nin}[0]{\noindent}

\newcommand{\sub}[0]{\subseteq}
\newcommand{\sm}[0]{\setminus}
\renewcommand{\dots}[0]{,\ldots,}

\newcommand{\ov}[0]{\overline}

\newcommand{\A}[0]{{\cal A}}
\newcommand{\B}[0]{{\cal B}}
\newcommand{\cee}[0]{{\cal C}}
\newcommand{\D}[0]{{\cal D}}
\newcommand{\eee}[0]{{\cal E}}
\newcommand{\f}[0]{{\cal F}}
\newcommand{\g}[0]{{\cal G}}
\newcommand{\h}[0]{{\cal H}}

\newcommand{\K}[0]{{\cal K}}

\newcommand{\pee}[0]{{\cal P}}
\newcommand{\Q}[0]{{\cal Q}}
\newcommand{\R}[0]{{\cal R}}
\newcommand{\sss}[0]{{\cal S}}
\newcommand{\T}[0]{{\cal T}}

\newcommand{\Z}[0]{{\cal Z}}

\newcommand{\ra}[0]{\rightarrow}
\newcommand{\Ra}[0]{\Rightarrow}

\newcommand{\ZZ}[0]{{\bf Z}}

\newcommand{\Nn}[0]{{\mathbb N}}

\newcommand{\PPP}[0]{\mathbb P}

\newcommand{\cc}[0]{\mbox{{\sf c}}}

\newcommand{\ddd}[0]{\gl}

\newcommand{\aaa}[0]{\ga}
\newcommand{\bbb}[0]{\gb}
\newcommand{\hhh}[0]{m}
\newcommand{\mm}{\mbox{{\sf m}}}
\newcommand{\qq}{\mbox{{\sf q}}}
\newcommand{\qqq}{q}
\newcommand{\ww}{\mbox{{\sf w}}}
\newcommand{\zz}[0]{\mbox{{\sf z}}}
\newcommand{\ghi}{\gc}

\newcommand{\0}[0]{\emptyset}

\renewcommand{\qed}[0]{\begin{flushright} \rule{2mm}{3mm} \end{flushright}}

\newcommand{\C}[2]{{{#1}\choose{{#2}}}}
\newcommand{\Cc}[0]{\tbinom}
\newcommand{\ga}[0]{\alpha }
\newcommand{\gb}[0]{\beta }
\newcommand{\gc}[0]{\gamma }
\newcommand{\gd}[0]{\delta }
\newcommand{\gD}[0]{\Delta }

\newcommand{\gl}[0]{\lambda }
\newcommand{\gL}[0]{\Lambda}
\newcommand{\go}[0]{\omega}
\newcommand{\gO}[0]{\Omega}

\newcommand{\gs}[0]{\sigma}

\newcommand{\gz}[0]{\zeta}
\newcommand{\eps}[0]{\varepsilon }
\newcommand{\vt}[0]{\vartheta}
\newcommand{\vs}[0]{\varsigma}
\newcommand{\vu}[0]{\upsilon}

\newcommand{\vp}[0]{\varphi}

\newcommand{\nex}[0]{_{\ov{x}}}

\newcommand{\sugg}[1]{}

\newcommand{\red}[1]{\textcolor{Red}{#1}}

\newcommand{\comments}[1]{}
\usepackage[usenames]{color}
\begin{document}

\renewcommand{\thefootnote}{\fnsymbol{footnote}}
\footnotetext{AMS 2010 subject classification:  05D40, 05D05, 05C65}
\footnotetext{Key words and phrases:  Erd\H{o}s-Ko-Rado property, random hypergraph}

\title{
On Erd\H{o}s-Ko-Rado for random hypergraphs I\footnotemark}
\author{A. Hamm and J. Kahn}
\date{}
\footnotetext{ $^*$Supported by NSF grant DMS1201337}

\maketitle

\begin{abstract}
A family of sets is {\em intersecting} if no two of its members are disjoint,
and has the {\em Erd\H{o}s-Ko-Rado property} (or {\em is EKR}) if each of
its largest intersecting subfamilies has nonempty intersection.

Denote by $\h_k(n,p)$ the random family
%on $[n]=\{1\dots n\}$
in which each $k$-subset of $\{1\dots n\}$ is present
with probability $p$, independent of other choices.
A  question first studied by
Balogh, Bohman and Mubayi asks:
\[
\mbox{{\em for what $p=p(n,k)$ is $\h_k(n,p)$ likely to be EKR?}}
\]
Here, for fixed $c<1/4$, and
$k< \sqrt{cn\log n}$
we give a precise answer to this question, characterizing those sequences
$p=p(n,k)$ for which
\[
\Pr(\h_k(n,p) ~\mbox{is EKR}) \ra 1 ~~\mbox{as} ~n\ra \infty.
\]

\end{abstract}

\section{Introduction}\label{Intro}

\bn
One of the most interesting combinatorial trends of the last couple decades
has been the investigation of ``sparse random" versions of
some of the
classical theorems of the subject---that is, of
the extent to which
these theorems hold in a random setting.
This issue has been the subject of some
spectacular successes, particularly those related to the
theorems of Ramsey \cite{Ramsey}, Tur\'an \cite{Turan} and Szemer\'edi \cite{Sz};
see \cite{FR86,BSS,RR, KLR}
for origins and, e.g., \cite{Conlon-Gowers, Schacht,DKTuran}
(or the survey \cite{Rodl-Schacht})
for a few of the most recent developments.

Here we are interested in what can be said in this vein
for the Erd\H{o}s-Ko-Rado
Theorem \cite{EKR},
another cornerstone of extremal combinatorics.
This natural question has already been considered by
Balogh, Bohman and Mubayi \cite{BBM}, and we first
quickly recall a few notions from that paper.

In what follows $k$ and $n$ are always positive integers with
$n>2k$.  As usual we write $[n]$ for $\{1\dots n\}$ and
$\C{V}{k}$ for the collection of $k$-subsets of the set $V$.
A $k$-{\em graph} (or $k$-{\em uniform hypergraph}) on $V$
is a subset (or multisubset), $\h$, of $\C{V}{k}$.
Members of $V$ and $\h$ are called {\em vertices} and {\em edges} respectively.
Here we will always take $V=[n]$ and write $\K$ for $\C{V}{k}$.
For a $k$-graph $\h$ on $V$ and $x\in V$ we use
$d_\h(x)$ for the {\em degree} of $x$ in $\h$ (the
number of edges of $\h$ containing $x$) and
$\gD_\h$ for the maximum of these degrees.
We also write $\h_x$ for the set of edges containing $x$,
called the {\em star} of $x$ in $\h$.

A collection of sets is {\em intersecting}, or a {\em clique},
if no two
of its members are disjoint.
The Erd\H{o}s-Ko-Rado Theorem says that for any $n$ and $k$
as above, the maximum size of an intersecting $k$-graph on $V$
is $\C{n-1}{k-1}$ and, moreover, this bound is achieved only
by the stars.
%{\em principal} $k$-graphs, those consisting of
%all $k$-sets containing some fixed element of $V$.

Following \cite{BBM}, we say that $\h\sub \K$
satisfies {\em (strong) EKR} if
every largest clique of $\h$ is a star;
thus the Erd\H{o}s-Ko-Rado  Theorem says $\K$ itself satisfies EKR.
(In \cite{BBM}
$\h$ is also said to satisfy
{\em weak EKR} if {\em some} largest clique is a star, but this slightly weaker
condition
will not concern us here.)

In what follows we use $\h=\h_k(n,p)$ for the random $k$-graph
on $V$ in which members of $\K$ are present independently,
each with probability $p$.
As suggested above, we are interested in understanding when
EKR holds for $\h$; a little more formally:
\begin{question}\label{Q1}
For what $p_0=p_0(n,k)$
is it true that $\h$ satisfies EKR a.s. provided $p> p_0$?
\end{question}
\nin
(As usual, an event---really a sequence of events parameterized by $n$ (say)---holds
{\em almost surely} (a.s.) if its probability tends to 1 as $n\ra\infty$.
Note that here we are thinking of $k$ as a function of $n$
({\em cf.} the paragraph following Theorem~\ref{MT}).)

Notice that EKR is not an increasing property (that is, it is not
preserved by addition of edges) and that, for given $n$ and $k$,
\beq{fnkp}
f_{n,k}(p):=\Pr(\h_k(n,p) ~\mbox{satisfies EKR})
\enq
is not increasing in $p$.
For instance, for sufficiently tiny $p$ (depending on $n$ and $k$) it will usually
be the case that {\em every} clique is contained in a star.
In view of this non-monotonicity,
it is natural to define a {\em threshold} for the property EKR to be
the least $p_0=p_0(n,k)$ satisfying
\beq{threshold}
f_{n,k}(p)\geq 1/2 ~~\forall p\geq p_0.
\enq
%(Alternatively, as in \cite{Erdos-Renyi} (or e.g. \cite[p. 18]{JLR}, we could say
%$p_0$ is {\em a} threshold for $Q$ if $f_{n,k}(p) > 1-o(1)$ for $p =\go(p_0$
%and this is not true if we replace $p_0$ by any $p_0'=o(p_0)$.  \red{[OY; skip?]}
%The present version
(This follows the usage in \cite{Kahn-Kalai} (e.g.), which takes
the ``threshold" for an {\em increasing} property $Q$ to be
the unique $p$ for which the ``$p$-measure" of
$Q$ is 1/2.)
\iffalse
In general---though we will sometimes do better---we
tend to regard determination of this threshold to within a constant factor as a satisfactory
answer to Question~\ref{Q1}.
\red{[But this may no longer be the best way to present, since we're asking for
a.s. behavior.]}
\fi

\medskip
For the most part we will not review the contents of \cite{BBM}.
The focus there
is mainly on small $k$;
roughly speaking, the authors give fairly complete results for
$k =o(n^{1/3})$ and more limited information for $k$ up to
$ n^{1/2-\eps}$ with $\eps>0$ fixed.
% (about which we will say a little more below).

The nature of the problem changes around
$k=n^{1/2}$, since for $k$ smaller than this, two random
$k$-sets are typically disjoint, while the opposite is
true for larger $k$.
Heuristically we may say that the problem becomes more interesting/challenging
as $k$ grows and the potential
violations of EKR proliferate (though increasing $k$ does narrow the range of $p$ for which
we {\em expect} EKR to hold).
At any rate, \cite{BBM} had (as noted there) little to say about $k$ larger than
$\sqrt{n}$ (or, indeed, $k>n^{1/2-\eps}$).
Here (in Theorem~\ref{MT'}) we precisely settle the problem for
$k $ up to and even a little beyond
$\sqrt{n}$.

As in \cite{BBM}, we will usually find it convenient to work, not directly with $p$, but
with $\vp :=p\C{n-1}{k-1}$, the expected degree of a vertex
(called $\rho$ in \cite{BBM}); this seems more natural as we are most interested
in situations where $p$ is tiny while the value of $\vp$ is more reasonable.
Throughout the paper we take $\mm=\E|\h|=\vp n/k$, $\gD=\gD_\h$ (the maximum degree in $\h$)
and
\beq{q}
\qq=\Pr(A\cap B\neq \0),
\enq
where $A$ and $B$ are chosen uniformly and independently from $\K$.
The next assertion will account for most of the work in our proof of
Theorem~\ref{MT'}.

\begin{thm}\label{MT}
For any fixed $c<1/4$, if
\beq{rangeofk}
k< \sqrt{cn\log n}
\enq
and $\vp$ is such that
\beq{avoid.generic}
\mbox{$\C{\mm}{\gD}\qq^{\C{\gD}{2}} <o(1) ~~$ a.s.,}
\enq
then $\h$ satisfies EKR a.s.
\end{thm}
\nin
(Recall $\C{a}{b} = (a)_b/b!:= a(a-1)\cdots (a-b+1)/b!$
for $a\in \mathbb R$ and $b\in \mathbb N$.)
Note that (here and usually in what follows) $n$ is a ``hidden parameter"; thus
in Theorem~\ref{MT}, $k$ and $\vp$ are
functions of $n$ and, for example, both ``$o(1)$" and ``a.s." in \eqref{avoid.generic} refer to $n\ra \infty$.
It may also be helpful to rephrase \eqref{avoid.generic}:
Given $\vp=\vp(n)$, set, for $t\in \Nn$,
\beq{f(t)}
\gL(t)=\gL_\vp(t)=\C{\mm}{t}\qq ^{\C{t}{2}}.
\enq
Then \eqref{avoid.generic} says:
\[
\mbox{$\exists$ $\eps=\eps(n)=o(1)$ such that $\Pr(\gL(\gD)>\eps)\ra 0$ as $n\ra\infty$.}
\]

Its meaning---the reason it is a natural assumption in Theorem~\ref{MT}
%of \eqref{avoid.generic}
---is as follows.
We think of $q^{\C{t}{2}}$ as the ideal
value of the probability that random (independent) $k$-sets $A_1\dots A_t$ form a clique (it
would be the actual value if the events $\{A_i\cap A_j\neq\0\}$
were independent).
Thus, since $|\h|$ is usually close to $\mm$,
the left side of \eqref{avoid.generic} may be thought of as
the expected number of ``generic" $\gD$-cliques in $\h$,
and we should perhaps not expect EKR to hold if this number is not small.

At least for $k$ as in \eqref{rangeofk}, this intuition turns out to be
correct; that is,  \eqref{avoid.generic}
is  essentially {\em necessary} for the conclusion of Theorem~\ref{MT}.
Here we should be a little careful:  since all cliques of size at
most 2 are trivial (that is, are contained in stars),
failure of \eqref{avoid.generic} with $\gD\leq 2$
does {\em not} suggest failure of EKR.
We accordingly define (again, given $\vp$)
\[  %\beq{gL'}
\gL'(t)=\gL'_\vp(t)=\left\{\begin{array}{ll}
0&\mbox{if $t\leq 2$,}\\
\gL(t)&\mbox{otherwise.}
\end{array}\right.
\]  %\enq
\begin{thm}\label{MT'}
For c and $k$ as in Theorem~\ref{MT} and any $\vp$ ($=\vp(n)$),
\beq{NandS}
\mbox{$\h$ satisfies EKR a.s. iff $\gL'(\gD) <o(1)$ a.s.}
\enq
\end{thm}
%\begin{thm}\label{MT'}
%For c and $k$ as in Theorem~\ref{MT} and any $\vp$ ($=\vp((n)$),
%$\h$ satisfies EKR a.s. iff
%\beq{NandS}
%\mbox{$\gL'(\gD) <o(1)$ a.s.}
%\enq
%\end{thm}
\nin
(That Theorem~\ref{MT} implies sufficiency of the condition in \eqref{NandS}
is easy but not quite tautological and will be discussed in Section~\ref{Nec}.)

\medskip
It is not hard to read off threshold information
from Theorem~\ref{MT'}
(with ``threshold" as in \eqref{threshold}, here translated
to the corresponding $\vp_0$); for example, for $k =\sqrt{\gz n}\gg \sqrt{n}$ (satisfying \eqref{rangeofk}),
we have $\vp_0\sim e^\gz\log n$.
Other special cases include the main positive results on EKR given in \cite{BBM}, those in
parts (i), (ii) and (iv) of their Theorem 1.1.
(We do use some of these in Section~\ref{Small}, but this could easily be avoided.)

Recent work of Balogh {\em et al.}
\cite{BDDLS}
provides results for $k$ up to $n/4$ but with nothing like the present accuracy.
(For $k$ as in \eqref{rangeofk}
%Theorem~\ref{MT} 
their upper bound on $\vp_0$ is of the form $e^{O(k)}\log n$.)
%exceeds the value established here by a factor $e^{\Theta(k)}$.)

\medskip
We believe Theorem~\ref{MT} is true with ``$c< 1/4$" replaced by ``$c<1/2$."
It is {\em not} true above this, roughly because:
for $k=\sqrt{cn\log n}$, with $c>1/2$, \eqref{avoid.generic} first occurs
at $\vp\approx\log n/\log(1/\qq )\sim n^c\log n$, where
%it will be the case that
(typically) all degrees are close to $\vp$ and for
each vertex $x$ the number of edges of $\h\sm \h_x$ meeting all edges
of $\h_x$ is about $\vp (n/k) \qq ^{\vp} \approx n^{c +1/2 -1}=n^{\gO(1)}$,
meaning stars are unlikely even to be maxi{\em mal} cliques.

This is, of course, reminiscent of the Hilton-Milner Theorem \cite{H-M}, which says that
the largest {\em nontrivial} cliques in $\K$ are those of the form
$\{A\}\cup \{B\in \K_x:B\cap A\neq\0\}$ (with $A\in \K$ and $x\in V\sm A$).
It seems not impossible that ``generic" and ``HM" cliques
are the main obstructions to EKR in general;
a precise, if optimistic, statement to this effect is:
\begin{conj}\label{EKRC}
If $k$ and $\vp$ are such that \eqref{NandS} holds and
\beq{noHM}
\mbox{$\h$ a.s. does not contain a Hilton-Milner family of size $\gD$,}
\enq
then $\h $ a.s. satisfies EKR.
\end{conj}
\nin
(Getting from this to the asymptotics of $p_0$ is routine.
Essentially---not quite literally---excluding HM families of size $\gD$ is promising
that stars are maxi{\em mal} cliques, and a slight weakening of
Conjecture~\ref{EKRC}, resulting in an unnoticeable change in the corresponding $p_0$, would replace
\eqref{noHM} by the assumption that this is true a.s.)

%%%%%%%%%%%%%%%%%%%%%%%%%%%%%%%%%%%%%%%%%%%%%%%%%%%%%%%%
%%%%%%%%%%%%%%%%%%%%%%%%%%%%%%%%%%%%%%%%%%%%%%%%%%%%%%%%

\iffalse

\begin{question}\label{Qwild}
Is it true that (for any k) if $\vp$ satisfies \eqref{avoid.generic} and
\beq{hxmaxl}
\mbox{a.s. every $\h_x$ is a maximal clique in $\h$,}
\enq
then $\h $ satisfies EKR a.s.?
\end{question}
\nin
(Understanding when \eqref{hxmaxl} holds is routine.)
%.

\fi

%%%%%%%%%%%%%%%%%%%%%%%%%%%%%%%%%%%%%%%%%%%%%%%%%%%%%%%%
%%%%%%%%%%%%%%%%%%%%%%%%%%%%%%%%%%%%%%%%%%%%%%%%%%%%%%%%

\medskip
In a companion paper \cite{HK2}, using methods completely different
from those employed here,
we jump to the other end of the spectrum, taking $k$ to be as
large as possible:

\begin{thm}\label{MT2}
There is a fixed $\eps>0$ such that if $n=2k+1$ and $p>1-\eps$, then
$\h $ satisfies EKR a.s.
\end{thm}

\nin
This was prompted by Question 1.4 of \cite{BBM}, {\em viz.}
\begin{question}\label{BBMQ}
Is it true that for $k\in (n/2-\sqrt{n},n/2)$ and $p=.99$,
EKR (or weak EKR) holds a.s. for $\h$?
\end{question}

\nin
Conjecture~\ref{EKRC} would say that
Theorem~\ref{MT2} remains true for $p$ at least about 3/4.
(Theorem~\ref{MT2} could presumably be extended to the full
range of $k$ covered by Question~\ref{BBMQ}, but this appears to be far short of the truth
if $n\geq 2k+2$, so seems of less interest.)

\bigskip
The rest of this paper is organized as follows.
The problem is most interesting when
\beq{kbig}
k> n^{1/2-o(1)}.
\enq
The bulk of our discussion of Theorem~\ref{MT} (Sections~\ref{Sketch} and \ref{PL1}-\ref{Large})
will deal exclusively with this range, while Section~\ref{Small} handles smaller $k$.
(Section~\ref{Prelim} reviews a few standard tools and Section~\ref{Gen} gives
some generalities that will apply to both regimes.)

In proving Theorem~\ref{MT} for $k$ as in \eqref{kbig} we will find it better
to deal first
with $\vp$ not too far above the ``threshold"---this will account for
most of our work---and then treat
larger $\vp$ mostly by a reduction to what we've established for smaller.
We thus begin in Section~\ref{Sketch} with an outline of the argument for small $\vp$,
in particular deriving Theorem~\ref{MT} in this range from three main assertions,
Lemmas~\ref{ClaimA}-\ref{ClaimC}.
These are proved in Sections~\ref{PL1}-\ref{PL3}
following the preparations of Sections~\ref{Prelim} and \ref{Gen}.
Section~\ref{Large} then gives the extension to large $\vp$ and, as noted above,
Section~\ref{Small} deals with small $k$.
Section~\ref{Nec} deals mainly with necessity of the condition in \eqref{NandS}.
This turns out to be interesting and considerably trickier than one might expect; still,
the paper being already too long, we will give the argument somewhat sketchily
and only for $k$ as in \eqref{kbig}.  (The problem gets easier as $k$ shrinks.)

\mn
{\bf Usage.}

As already mentioned we take
%Throughout the paper we take
$V=[n]$, $\K=\C{[n]}{k}$ and $\h=\h_k(n,p)$.
In addition we set $M=\C{n-1}{k-1}$ (so $\vp = Mp$) and $m=|\h|$
(a random variable with mean $\mm$).
We use $v,w,x,y,z$ for members of $V$.
For a hypergraph $\g$,
we let $\g\nex=\g\sm\g_x$
(recall $\g_x=\{A\in \g:x\in A\}$).

We use $d_\g(x)$ for the degree of $x$ in $\g$, and similarly for $d_\g(x,y)$,
and, where not otherwise specified, take $d(\cdot)$ to mean $d_\h(\cdot)$.
(As already stated, we use $\gD$ for $\gD_\h$.)

We write $B(\ell,\ga)$ for a random variable
with the binomial distribution ${\rm Bin}(\ell,\ga)$,
$\log $ for $\ln$ and $\C{a}{\leq b} $ for $\sum_{i\leq b}\C{a}{i}$.
We use standard asymptotic notation (``big Oh" etc.), but will also sometimes use
$a\asymp b$ for $a=\Theta(b)$ and $a\ll b$ for $a =o(b)$.
We assume throughout that $n$ is large enough to support our arguments and,
following a standard abuse, usually pretend large numbers are integers.

\section{Main points}\label{Sketch}

From now until the end of Section~\ref{Large}
we fix $c= 1/4-\eps$ in Theorem~\ref{MT}.
Also, as noted above, the present section assumes $k$ satisfies \eqref{kbig} (as well as \eqref{rangeofk}).

As stated earlier, most of our work will deal with $\vp$ fairly near the ``threshold."  Though the problem should
become easier as $\vp$ grows, some parts of the main argument below break down
for larger $\vp$; this could perhaps be remedied, but we have found it easier to
first deal directly with smaller $\vp$ and then use what we've learned to handle larger values.
%handle large $\vp$ by comparison with smaller values.
%
(A disadvantage of this approach is that it necessarily gives much weaker bounds on the probability
that EKR fails than one might {\em hope} to establish using a more direct argument.)

We thus begin in this section with an outline of where we are headed
in the ``small $\vp$" regime.
As we will see, the ``threshold" $\vp_0$ ($:=Mp_0$)
is around $\log n/\log(1/\qq)$, and as a cutoff for ``small" we set
(not a delicate choice)
\beq{vp*}
\vp^* = \frac{\log^3n}{\log(1/\qq)}.
\enq

We assume in this section (and again in parts of Section~\ref{Gen} and all of Sections~\ref{PL1}-\ref{PL3}) that
$\vp\leq \vp^*$
(a restriction which could be relaxed considerably without invalidating the present argument).
Thus we want to show
\beq{smallvp}
\mbox{for $\vp\leq \vp^*$ satisfying \eqref{avoid.generic}, $\h$ satisfies EKR a.s.}
\enq
(It {\em is} true that in this regime the problem is most delicate when $\vp$ is more or less
at the ``threshold"; in particular it is only here---see the proof of Lemma~\ref{ClaimC}---that
we must make precise use of \eqref{avoid.generic}.)

\medskip
Call a clique {\em trivial} if it is contained in a star.
We will show below that
there are integers $\aaa=\aaa(n,\vp)\leq \bbb=\bbb(n,\vp)$
satisfying, {\em inter alia},
\beq{gDinab}
\gD\in[\aaa,\bbb] ~~\mbox{a.s.}
\enq
and
\beq{chooset''}
\gL(\aaa) = o(1).
\enq

\mn
Thus Theorem~\ref{MT} would follow if we could show that $\h$
a.s. does not contain a nontrivial clique of size $\aaa$,
but this is not quite true; for example,
if $d(x)=\gD$ is significantly larger than $\aaa$---say
closer to $\bbb$ than $\aaa$---then an $A\in \h\sm\h_x$
{\em typically} misses fewer than $\gD-\aaa$ edges of $\h_x$,
in which case $\{A\}\cup \{B\in \h_x:B\cap A\neq\0\}$ is a nontrivial
clique of size greater than $\aaa$.

A natural way to address this is to compare each clique possessing a
sufficiently high degree vertex, say $x$, directly with the star $\h_x$.
This idea is implemented in the first of the following three lemmas;
these assertions will easily yield \eqref{smallvp} and will also do most of the work
when we come to larger $\vp$.
(To be clear, the lemmas will depend on further properties of $\aaa$ and $\bbb$
to be established below.)

Set
\beq{ghi}
\ghi=\min\{\aaa,\vp^*/3\},
\enq
\beq{s0}
\tau = (1-\eps)\ghi
\enq
and
\beq{gl}
\gl = \max\left\{\frac{\sqrt{\log n}}{\log (1/\qq )}, 2\sqrt{\frac{\log n}{\log (1/\qq )}}~\right\}.
\enq
%\beq{gl}
%\gl = \max\{\sqrt{\log n}/\log (1/\qq ), 2\sqrt{\log n/\log (1/\qq )}\}.
%\enq
(The actual values are not needed in this section.
One should think of $\ghi=\aaa$; the technical $\vp^*/3$ will be needed
for the reduction in Section~\ref{Large}.)

%
%Note in particular that, in view of \eqref{chooset''}, we have
%\beq{lamgam}
%\gL(\ghi)=o(1).
%\enq

\begin{lemma}\label{ClaimA}
A.s. there do not exist (in $\h$) a nontrivial clique $\cee$ and vertex x such that
$|\cee|\geq d(x)$,
$d_\cee(x)\geq \tau $,  and
either $|\cee|\geq \aaa$ or $|\cee\nex|\geq 2/\eps$.
\end{lemma}

\begin{lemma}\label{ClaimB}
A.s. $\h$ does not contain a clique with two vertices
of degree at least $\gl$.
\end{lemma}

\begin{lemma}\label{ClaimC}
A.s.
$\h$
does not contain a clique of size $\ghi$ with
at most one vertex of degree greater than $\gl$ and maximum degree less than $\tau $.
\end{lemma}

\nin
(For perspective we remark that
Lemmas~\ref{ClaimA} and \ref{ClaimC} are the main points;
Lemma~\ref{ClaimB} just makes our lives a little easier
when we come to Lemma~\ref{ClaimC}.)

\medskip
Lemmas~\ref{ClaimA}-\ref{ClaimC}
easily imply \eqref{smallvp},
as follows.
Since $\Pr(\gD<\aaa) =o(1)$ (see \eqref{gDinab}), it is enough to show that $\h$
a.s. does not contain a nontrivial clique $\cee$ with $|\cee|\geq \gD\geq \aaa$.
But if $\gD\geq \aaa$ and $\h$ does contain such a $\cee$, then at least
one of the following occurs.

\mn
(a)  There is an $x$ with $d_\cee(x)\geq \tau $
(and $|\cee|\geq \gD\geq\max\{\aaa,d(x)\}$), so $x,\cee$ are as in Lemma~\ref{ClaimA}.

\mn
(b)
There are two vertices with degree at least $\gl$ in $\cee$.

\mn
(c)  There is at most one vertex $x$ with $d_\cee(x)\geq \gl$
and none with $d_\cee(x)\geq \tau $, so (since $\aaa\geq \ghi$) $\cee$ is as in Lemma~\ref{ClaimC}.

\mn
But according to Lemmas~\ref{ClaimA}-\ref{ClaimC}, each of (a)-(c) occurs with probability $o(1)$,
so we have \eqref{smallvp}.\qed

As the reader may have noticed, this derivation would remain valid if we
dropped the alternative ``$|\cee\nex|\geq 2/\eps$" in Lemma~\ref{ClaimA} and replaced $\gc$
by $\aaa$ in Lemma~\ref{ClaimC};
the stated versions of these lemmas will be needed for dealing with larger $\vp$ in Section~\ref{Large}.

\medskip
We now pause to fill in some preliminaries.

\section{Negative association and large deviations}\label{Prelim}

Some parts of the analysis below seem most conveniently handled using
the notion of
negative association, regarding which we just
recall what little we need, in particular confining ourselves to
$\{0,1\}$-valued r.v.'s; see e.g. \cite{Pem,DR} for further background.

Recall that events $\A,\B$ in a probability space
are \emph{negatively
correlated}
(denoted $\A\downarrow \B$)
if $\Pr(\A\B)\leq\Pr(\A)\Pr(\B)$.
Given a set $S$, set $\gO=\gO_S=\{0,1\}^S$ and
recall that $\A\sub \Omega$
is {\em increasing}
if $x\geq y\in\A~\Ra x\in\A$ (where ``$\geq$" is product order on $\Omega$).
Say $i\in S$ {\em affects} $\A\sub \Omega$ if there are $\eta\in \A$
and $\nu\in\Omega\sm\A$ with $\eta_j=\nu_j ~\forall j\neq i$,
and write $\A\perp \B$ if no $i\in S$ affects both $\A$ and $\B$.

Now suppose $(X_i:i\in S)$ is drawn from some probability distribution on $\gO$.
%and write $\mu$ for the law of $X$ (thus $\mu(\A)=\Pr(X\in \A)$ for $\A\sub \gO$).
The $X_i$'s are said to be {\em negatively associated}
(NA) if
$\A\downarrow \B$
%$\mu(\A\B)\leq \mu(\A)\mu(\B)$
whenever $\A,\B$ are increasing and $\A\perp \B$.
If $Q_i$ are events whose indicators are NA then we also say
that the $Q_i$'s themselves are NA.

The following observation is surely not news,
but as we don't know a reference we give the easy proof.

\begin{prop}\label{NAProp}
Suppose that for some $V_1\dots V_s\sub V$ and $\ell_1\dots \ell_s$,
$A_1\dots A_s$ are chosen independently with $A_j$ uniform from
$\C{V_j}{\ell_j}$.
Then the r.v.'s $X_{vj}= {\bf 1}_{\{v\in A_j\}}$
($v\in V$, $j\in [s]$) are negatively associated.
\end{prop}

\mn
{\em Proof.}  ({\em Cf.} \cite[Prop. 12]{DR}.)
For each $j$ the vector $(X_{vj}:v\in V)$ is chosen uniformly from
the strings of weight $\ell_j$ in $\{0,1\}^{V_j}$, implying that the r.v.'s $X_{vj}$ ($v\in V$) are
NA.  (This is standard and easy, though we couldn't find it in writing.
A stronger and far more interesting statement is the main result of \cite{vdBJ}.)
%
%(This is an instance of, for example, \cite[Theorem 1.3]{Pem},
%but is surely older---probably folkloric---and is easily proved directly.)
%
We may thus apply \cite[Proposition 8]{DR}, which says that if the collections
$\{X_{vj}:v\in V\}$ ($j\in [s]$) are mutually independent and each is NA,
then the entire collection $\{X_{vj}\}$ is also NA.\qed

We will use Proposition~\ref{NAProp} in conjunction with the following trivial observations.

\begin{prop}\label{NAProp'}
If the r.v.'s $X_1\dots X_m$ are NA,
$I_1\dots I_r$ are disjoint subsets of $[m]$, and $Q_j$ is an increasing event
determined by $\{X_i:i\in  I_j\}$, then $Q_1\dots Q_r$ are NA.
\end{prop}

\begin{prop}\label{NAProp''}
If the events $Q_i$ are NA, then $\Pr(\cap Q_i)\leq \prod \Pr(Q_i)$.
\end{prop}

One virtue of negative association lies in the fact that
``Chernoff-type" large deviation bounds for random variables $X=\sum X_i$, where
$X_1,\ldots$ are
independent Bernoullis, remain valid under the (weaker) assumption
that the $X_i$'s are
negatively associated.
As far as we know, this was first observed by Dubhashi and Ranjan \cite[Proposition 7]{DR};
it is gotten via the usual argument (Markov's inequality applied to $\exp[tX]$;
see e.g. \cite[pp. 26-28]{JLR}),
with the identity $\E e^{tX} =\prod \E e^{tX_i}$ replaced by the inequality
$\E e^{tX} \leq \prod \E e^{tX_i}$.
In particular this gives the following bounds
(see for example \cite[Theorem 2.1 and Corollary 2.4]{JLR}).
\begin{thm}\label{Chern}
Suppose $X_1\dots X_m$ are negatively associated ${\rm Ber}(p)$ r.v.s,
$X=\sum X_i$ and $\mu=\E X$.  Then for any $\gl\geq 0$,
\begin{eqnarray}
\Pr(X >\mu+\gl)& < &\exp [- \tfrac{\gl^2}{2(\mu+\gl/3)}],\label{Chern1}\\
\Pr(X < \mu-\gl) &<& \exp [- \tfrac{\gl^2}{2\mu}],\nonumber
\end{eqnarray}
and for any $K>1$,
\beq{Chern3}
\Pr(X > K\mu)< [e^{K-1}K^{-K}]^\mu.
\enq
\end{thm}

\begin{cor}\label{CChern}
The inequality
\eqref{Chern3} still holds if instead of $\E X=\mu$ (in  Theorem~\ref{Chern})
we assume only $\varrho:=\E X\leq \mu$.
\end{cor}

\nin
{\em Proof.}
We have (using \eqref{Chern3} for the inequality)
\begin{eqnarray*}
\Pr(X > K\mu) &=& \Pr(X> (K\mu/\varrho)\varrho)\\
&< &[e^{K\mu/\varrho-1}(K\mu/\varrho)^{-K\mu/\varrho}]^\varrho
=
e^{K\mu-\varrho}K^{-K\mu}(\mu/\varrho)^{-K\mu}.
\end{eqnarray*}
The last expression is equal to the bound in \eqref{Chern3} when $\mu=\varrho$
and is easily seen to be decreasing in $\mu\geq \varrho$ (provided $K\geq 1$).\qed

\section{Generics}\label{Gen}

This section establishes basic properties
of some of the
parameters we will be dealing with, in particular showing that
$\h$ a.s. satisfies a few general properties
whose failure can then be more or less ignored in what follows.

\medskip
To begin, we should say something about the intersection probability $\qq$ (defined in \eqref{q}).
We have $\qq=1-\vt$ with
\beq{qvt}
\vt= \tfrac{(n-k)_k}{(n)_k}\sim e^{-k^2/n},
\enq
where, as usual, $(b)_a=b(b-1)\cdots (b-a+1)$.
(The ``$\sim$" is valid provided $k=o(n^{2/3})$.)
This gives the asymptotics of $\qq $ for $k=\gO(\sqrt{n})$;
in particular for $k\gg\sqrt{n}$ we have
\beq{log1q}
\log (1/\qq )\sim e^{-k^2/n}.
\enq
For $k \ll \sqrt{n}$ we instead have
\beq{qsmallk}
\qq \sim k^2/n
\enq
(since, with $X_v={\bf 1}_{\{v\in A\cap B\}}$,
\[
k^2/n=\sum\E X_v\geq \qq \geq \sum\E X_v-\sum\E X_vX_w > k^2/n-\Cc{n}{2}(k/n)^4).
\]
Note that in any case we have
\beq{vp*'}
\vp^*< n^{1/4-\eps+o(1)}.
\enq

\medskip
We will usually be dealing with situations in which $\qq$ is
slightly perturbed by information on how relevant $k$-sets meet
%how relevant members of $\K$ meet
some small subset of $V$.  This negligible effect is handled by the next observation.

\begin{prop}\label{q0prop}
Fix $W\sub V$ of size at most $\ww\ll n/\log n$ and $B\in \C{V}{k}$, and let $A$ be
uniform from $\C{V}{k}$.  Then conditioned on any value of $A\cap W$,
\[
\Pr(A\cap (B\sm W)\neq \0)<
 (1+2k^2\ww/(\qq n^2))\qq.
 \]
\end{prop}

\mn
{\em Proof.}
The probability is largest when $|W|=\ww$  and $B\cap W = A\cap W=\0$,
in which case its value is $\qqq =1-\vs$, with
$
\vs= \tfrac{(n-\ww-k)_k}{(n-\ww)_k}.
$
We have
\begin{eqnarray*}
\frac{\vt}{\vs} &=&\frac{(n-k)_k(n-\ww)_k}{(n)_k(n-\ww-k)_k}\\
&=&\prod_{i=0}^{k-1}\left(1+\frac{k\ww}{(n-i)(n-\ww-k-i)}\right)
= 1 + (1+o(1))\frac{k^2\ww}{n^2}~;
\end{eqnarray*}
that is, $\vt/\vs-1\sim k^2\ww/n^2$ ($=o(1)$ because of the bound on $\ww$).  Thus
\begin{eqnarray*}
\frac{\qqq }{\qq}-1 &=& \frac{\vt-\vs}{1-\vt} =\frac{1}{1-\vt}\left(\frac{\vt}{\vs}-1\right)\vs\\
&\sim & \frac{k^2\ww\vs}{(1-\vt)n^2}
\sim \frac{k^2\ww}{(1-\vt)n^2}e^{-k^2/n}.
\end{eqnarray*}
The lemma follows.
\qed

\medskip
In all that follows we assume $\vp$ satisfies \eqref{avoid.generic}.
At some (indicated) points in this section,
and again throughout Sections~\ref{PL1}-\ref{PL3}, we will also stipulate
that  $\vp\leq \vp^*$
(defined in \eqref{vp*}).

We also assume from now on that $\mm=\go(1)$, since cases with $\mm=O(1)$ are trivial:
if $\mm=o(1)$ then \eqref{avoid.generic} fails
(the l.h.s. is a.s. 1; actually in this case $\h$ is a.s. empty and does satisfy EKR),
while if $\mm =\Theta(1)$ then with probability $\gO(1)$ we have $\gD=|\h|=1$ and the l.h.s.
of \eqref{avoid.generic} is $\mm$ (so \eqref{avoid.generic} does not hold).
Recall $m=|\h|$.
Let $\psi=\psi(n)$ be some slowly growing function of $n$ (say $\psi = \log n$).
Theorem~\ref{Chern} (for independent Bernoullis)
says that a.s.
\beq{mm}
\hhh\in (\mm -\psi\sqrt{\mm},\mm +\psi\sqrt{\mm}).
\enq
We henceforth write $m_0 $ for $ \mm + \psi\sqrt{\mm}$.

It will often be convenient to replace $\h$ by
$A_1\dots A_m$ chosen {\em independently} from $\K$, a change which
makes little difference when $\mm$ is small:
\begin{prop}\label{indeptAs} If $m_0\ll \sqrt{|\K|}$
then for any property $\A$ and any c,
if,
\[\max_{m\models \eqref{mm}} |\Pr(A_1\dots A_m\models \A)-c| =o(1)
\]
(where the $A_i$'s are chosen uniformly and independently from $\K$), then
\[
\Pr(\h\models\A)=c+o(1).
\]
\end{prop}
\nin
(We will eventually need something more careful in a similar vein; see the paragraph preceding
Lemma~\ref{ClaimC'}.)

\mn
{\em Proof.}
For any $l$
the law of $\{A_1\dots A_l\}$ given that the $A_i$'s are distinct
is the same as that of $\h$ given $m=l$.
We may thus couple $\h$ and $\{A_1\dots A_m\}$ so that they coincide whenever the
$A_i$'s are distinct.  But the probability that they are {\em not}
distinct is at most
\[
\Pr(m\not\models \eqref{mm}) + \Pr(\mbox{$A_1\dots A_{m_0}$ are not distinct})
< o(1)+m_0^2/|K| =o(1)
\]
and the proposition follows.\qed

\medskip
From now until the ``coda" at the end of this section we assume that
\beq{vpm}
\vp > n^{-o(1)}.
\enq
As we will see in the coda, this is implied by \eqref{avoid.generic} if we assume \eqref{kbig}.
(Note that here we do {\em not} assume \eqref{kbig}, since we will also need parts of what follows
in Section~\ref{Small}, where \eqref{kbig} does not hold.)

\medskip
We next need to say something about the behavior of $\vp$ and $\gD$.
Recall that our default for degrees is $\h$;
thus, in addition to $\gD=\gD_\h$, we take $d_x=d(x)=d_\h(x)$ and $d(x,y)=d_\h(x,y)$.
The properties we need will be given in Proposition~\ref{things}
once we have introduced the parameters $\aaa$ and $\bbb$ mentioned earlier.

Let $\aaa_1$ and $\bbb$ be, respectively,
the largest integer with $\Pr(d_v\geq \aaa_1)\geq \psi/n$
and the smallest integer with $\Pr(d_v> \bbb)<1/(n\psi)$
(for any $v$).

Notice that $\gL(0)=1$ and (since $\gL(t)/\gL(t-1) = ((\mm-t+1)/t)\qq^{t-1}$ is decreasing in $t$)
there is some $t_0$ such that $\gL(t)$ is increasing up to $t_0$
and decreasing thereafter.
Thus \eqref{avoid.generic} says that there are
$\vs=\vs(n)$ and $\vu=\vu(n)$, both $o(1)$, such that $\Pr(\gL(\gD)>\vs) < \vu$.
Set $\aaa_2  :=\min\{t:\gL(t)\leq \vs\}$ and $\aaa =\max\{\aaa_1,\aaa_2\}$.

The promised Proposition~\ref{things} now collects properties of these
parameters that we will use repeatedly
in what follows, often without explicit mention.
\begin{prop}\label{things}
For $\aaa$, $\bbb$ as above:
\beq{thinga}
\aaa\leq\bbb;
\enq
\beq{thingb}
\gL(\aaa) =o(1);
\enq\beq{thingc}
\mbox{$\gD\leq \bbb$ a.s.; if $\vp \leq \vp^*$ then $\gD\geq \aaa$ a.s.};
\enq\beq{thingd}
\bbb/\vp < n^{o(1)};
\enq
\beq{thinge}
\aaa > (1-o(1))\log n/\log(1/\qq);
\enq\beq{thingf}
\mbox{if $\vp\leq\vp^*$ then $\bbb< (1+o(1))\vp^*$ ($< n^{1/4-\eps+o(1)}$).}
\enq
\end{prop}

\mn
(It is not hard to see that in fact $\aaa\sim\bbb$ in all cases
and $\bbb\sim \vp$ iff $\vp\gg \log n$.
What we actually use for the second part of \eqref{thingc} is $\aaa_1 k/n\ll 1$.)

\medskip
For the rest of the paper we set
$
\pee =\{\mbox{$m$ satisfies \eqref{mm}}\}\wedge
\{\gD\leq\bbb\},
$
noting that \eqref{thingc} and our earlier observation that
\eqref{mm} holds a.s. give
\beq{pee}
\Pr(\pee) = 1-o(1).
\enq

\mn
{\em Proof of Proposition~\ref{things}}.
The first assertion in \eqref{thingc} is immediate from the definition of $\bbb$.
From the definition of $\aaa_2$ we have
$  %\beq{chooset'}
\gL(\aaa_2) = o(1)
$  %\enq
(namely $\gL(\aaa_2)\leq \vs$) and $\gD\geq \aaa_2$ a.s.
(since $\Pr(\gD <\aaa_2)=\Pr(\gL(\gD)>\vs)<\vu$),
implying $\aaa_2\leq \bbb$.
This gives \eqref{thinga} and \eqref{thingb}.

\medskip
Let
$\bbb^* = \lceil\vp +\eta\rceil$, with
$\eta $ the positive root of
$x =\sqrt{2(\vp +x/3)(\log n+\log \psi)}$.
Then Theorem~\ref{Chern} gives (for any $v$)
\beq{beta}
\Pr(d_v > \bbb^*) < \exp[-\eta^2/(2(\vp +\eta/3))] = (n\psi)^{-1},
\enq
whence $\bbb\leq \bbb^*$.
(The bound is very crude for smaller values of $\vp$, but we have
lots of room in such cases.)
In particular, since $\eta = O(\max\{\sqrt{\vp \log n}, \log n\})$,
\eqref{vpm} now implies both \eqref{thingd} and \eqref{thingf}
(and $\bbb\sim \vp$ if $\vp\gg \log n$, but we don't need this).

For \eqref{thinge} we have
\begin{eqnarray*}
\gL(\aaa_2) &>& \exp[\aaa_2(\log (\mm/\aaa_2)-\tfrac{\aaa_2-1}{2}\log (1/\qq ))]\\
&>& \exp\left[\tfrac{\aaa_2}{2}((1-o(1))\log n  -\aaa_2\log (1/\qq ))\right]
\end{eqnarray*}
(since $\log (\mm/\aaa_2) > (1/2-o(1))\log n$, as
follows from $\mm = \vp n/k$, $\aaa_2\leq \bbb$ and \eqref{thingd}), and
combining this with \eqref{thingb} gives
$  %\beq{abig}
\aaa_2>(1-o(1))\log n/\log(1/\qq ).
$  %\enq

Finally, the second assertion in \eqref{thingc} is given by the following more general
statement, which we will need again in Section~\ref{Small}.
\begin{prop}\label{DeltaProp}
For any $n,k,\vp$ ($=Mp$) and $\theta\in \Nn$ satisfying $p=o(1)$ and $\theta=o(M)$:
if $\Pr(d_v\geq \theta) =\go(1/n)$ and $\theta k/n=o(1)$ then $\gD\geq \theta$ a.s.
\end{prop}
\nin
(The assumption $\theta =o(M)$ is a little silly:  for $k\geq 3$ it follows from $\theta k/n=o(1)$.
For \eqref{thingc}---note we already know $\gD\geq \aaa_2$ a.s.---the hypothesis $\aaa_1 k/n=o(1)$ follows from
$\vp k/n < n^{-1/4}$ and $\aaa_1/\vp\leq \bbb /\vp< n^{o(1)}$;
see \eqref{vp*'} and \eqref{thingd}.
For $k< n^{1/2-\gO(1)}$ and a {\em fixed} $\theta$, Proposition~\ref{DeltaProp} is \cite[Lemma 3.6]{BBM}.)

\mn
{\em Proof of Proposition~\ref{DeltaProp}.}
Let
$X_v={\bf 1}_{\{d_v\geq\theta\}}$ and $X=\sum X_v$.
We are assuming $\E X =\go(1)$, so to finish via the second moment method just need
\beq{ExvXw}
\E X_vX_w\sim \E^2X_v
\enq
(for $v\neq w$).  Letting $Z=d(v,w)$ we have
\beq{ExvXw'}
\mbox{$\E X_vX_w<\sum_{l\geq 0}\Pr(Z=l) \Pr^2(d_v\geq \theta-l).$}
\enq
(For equality we would replace $d_v$ by
$d(v,\ov{w}):= |\h_v\sm \h_w|$.)

Now, $Z$ is binomial with $\E Z< \vp k/n $, so
\beq{Z}
\Pr(Z=l)~ ~
(\leq \Pr(Z\geq l)) ~
< (\vp k/n)^l.
\enq

\mn
%(The ``$\Pr(Z\geq l)$" is for future reference.)
On the other hand, since $d_v\sim {\rm Bin}(M,p)$, we have,
for each $t\leq \theta$,
\beq{tvp}
\frac{\Pr(d_v=t-1)}{\Pr(d_v=t)} ~=~ \frac{t(1-p)}{(M-t+1)p} \sim t/\vp,
\enq
implying $\Pr(d_v\geq t-1)< (1+\theta/\vp)\Pr(d_v\geq t)$.
Thus (since $\theta k/n=o(1)$)
the sum in \eqref{ExvXw'} is asymptotic to its zeroth term, and we have \eqref{ExvXw}.

(We pickily add---to make sure that $\vp k/n=o(1)$---that we may assume $\theta\geq \vp$:
there is nothing to prove if $\theta =0$, and $\gD\geq \vp$ is easy if $\vp\geq 1$
(and $k=o(n)$, which follows from $\theta>0$ and $\theta k/n=o(1)$).)
\qed

We will also eventually (in Section~\ref{PL3}) need the easy:  if $\vp\leq \vp^*$ then
\beq{m0m}
\Cc{m_0}{\aaa} \sim \Cc{\mm}{\aaa}.
\enq
(The ratio of the left- and right-hand sides is
\[
\tfrac{(m_0)_\aaa}{(\mm)_\aaa} < (\tfrac{m_0-\aaa+1}{\mm-\aaa+1})^\aaa
< \exp[O(\psi\aaa/\sqrt{\mm})]
\]
and
$\psi\aaa/\sqrt{\mm}\leq \psi\bbb/\sqrt{\mm}< n^{-\eps+o(1)}$
(using $\mm =\vp n/k$, \eqref{vpm} and \eqref{thingf}).)

\medskip
For $x\in V$, let $W_x=\{y:d(x,y)\geq 2\}$ (a random set determined by $\h_x$).
Let $\R$ be the intersection of $\pee$ and the events $\{\gD\geq\aaa\}$,
\beq{dhxy}
\{d(x,y) \leq 8  ~~\forall x,y\}
\enq
and
\beq{Wx}
\{|W_x| < \max\{\vp^2k^2/n, 6\log n\} ~~\forall x\}.
\enq
Though defined here in general, $\R$ is only of interest when
$\vp$ is small:

\begin{prop}\label{HProp}
If $\vp\leq \vp^*$, then $\Pr(\R)=1-o(1)$.
\end{prop}

\nin
{\em Proof.}
We have already seen (in \eqref{pee} and \eqref{thingc}) that $\pee$ and
$\{\gD\geq \aaa\}$ hold a.s.
That \eqref{dhxy} does as well
follows ({\em via} the union bound) from
the fact, already noted in \eqref{Z}, that
$\Pr(d(x,y)\geq l)\ll n^{-l/4}$.
% (for any $\l \geq 0$).
%
To deal with \eqref{Wx}, it is, according to Proposition~\ref{indeptAs},
enough to show

\mn
{\em Claim.}  If $m$ satisfies \eqref{mm} and $A_1\dots A_m$ are chosen independently
(and uniformly)
from $\K$, then \eqref{Wx} holds a.s.

\mn
(where of course $d$ in the definition of $W_x$ now refers to $\{A_1\dots A_m\}$ rather than $\h$).

\mn
{\em Proof of Claim.}
For a given $x$ we have, for each $y\neq x$,
$\Pr(y\in W_x) < \C{m}{2}(k/n)^4< (1/2+o(1))(\vp k/n)^2$
(using $m\sim \mm =\vp n/k$), implying
$\E|W_x| < (1+o(1))\vp^2 k^2/(2n)$.
On the other hand,
the events $\{y\in W_x\}$ are NA
(by Propositions~\ref{NAProp} and \ref{NAProp'}), and a little calculation, with Corollary~\ref{CChern},
bounds the probability that a particular $x$ violates \eqref{Wx} by
$o(1/n)$.
(In more detail:
if $\mu :=\vp^2k^2/(2n) \geq 3\log n$, then
\eqref{Chern1} bounds the probability by $\exp[-(9/8)\log n]$;
otherwise $K:= 6\log n/\mu> 2$, and \eqref{Chern3} bounds the probability by
$(e^{K-1}K^{-K})^\mu = (e^{1-1/K}K^{-1})^{K\mu} \leq (\sqrt{e}/2)^{6\log n}=o(1/n)$.)

\qed

%\nin
%{\em Remark.}
%Proposition~\ref{WZprop}(a) gives a much
%stronger probability bound for an assertion similar to the above Claim.

\mn
{\bf Coda.}  Finally, we say why
the combination of \eqref{avoid.generic}
and \eqref{kbig} implies \eqref{vpm}.  Suppose instead that the first two
conditions hold but $\vp< n^{-\gO(1)}$.
Then $\gD <O(1)$ a.s.
But if $\gD=O(1)$, then $\qq> n^{-o(1)}$
(see \eqref{qsmallk}) implies
$\gL(\gD)  = \gO(\mm^\gD) n^{-o(1)}$, so that
\eqref{avoid.generic} implies $\mm < n^{o(1)}$
(note $\gD\geq 1$ a.s. since we assume
$\mm =\go(1)$).
But then (since $\mm =\vp n/k$ and we assume \eqref{kbig})
$\vp < n^{-1/2+o(1)}$, implying that in fact $\gD\leq 2$ a.s.

Now suppose $\gL(2) = o(1)$.  Then $k\ll \sqrt{n}$ (otherwise $\qq=\gO(1)$ and $\mm=o(1)$,
contrary to assumption), and $\gL(2)\asymp (\vp n/k)^2(k^2/n) =\vp^2 n$, implying
$\vp \ll n^{-1/2}$ and $\gD=1$ a.s.  But $\gL(1)=\mm$, so we contradict
\eqref{avoid.generic} .\qed

\section{Proof of Lemma \ref{ClaimA}}\label{PL1}

Here and in Sections~\ref{PL2} and \ref{PL3} we
take
\beq{ww}
\ww=\max\{\vp^2k^2/n, 6\log n\}
\enq
and
\beq{qqq}
\qqq =
(1+2k^2\ww/(\qq n^2))\qq;
\enq
thus $\ww$ is the bound on the $|W_x|$'s in \eqref{Wx}
(we will use it to bound a related quantity in Section~\ref{PL3})
and
$\qqq$ is the probability bound in Proposition~\ref{q0prop}.
We will need to say that $\qqq$ is close to $\qq$; here and in Section~\ref{PL2}
we could get by with, for example, $\log(1/\qqq )\sim \log(1/\qq)$,
but for the more delicate situation in Section~\ref{PL3} will need
\beq{qqqqq}
\qqq^{\C{\aaa}{2}}\sim \qq^{\C{\aaa}{2}}
\enq
(that is, $k^2\ww \aaa^2/(\qq n^2)=o(1)$; in fact,
$k^2\ww \aaa^2/(\qq n^2)< n^{-4\eps +o(1)}$ since
$\aaa < n^{1/4-\eps+o(1)}$ (see \eqref{thingf}), $\ww < n^{1/2-2\eps+o(1)}$
(see \eqref{vp*'}) and
$k^2/(\qq n) < 1+o(1)$.)

We will use (a) of the following observation in the present section and the variant (b) in Section~\ref{PL2}.

\begin{prop}\label{ABProp}
{\rm (a)}  Suppose $\A=\{A_1\dots A_d\}\sub \K_x$
satisfies
\beq{dandW}
\mbox{$d_\A(z)\leq 8 ~\forall z\in V\sm \{x\}~$ and $~|\{z\in V\sm \{x\}:d_\A(z)\geq 2\}|<\ww$.}
\enq
Then for $B$ uniform from $\K\nex$,
\[
\Pr(B\cap A_i\neq\0 ~\forall i\in [d]) < (1+o(1))\qqq ^d.
\]
{\rm (b)}  The same conclusion holds if
$\A\sub \{A\in \K_x :y\not\in A\}$ satisfies \eqref{dandW} and
$B$ is uniform from $\{A\in \K_y :x\not\in A\}$.
\end{prop}
\nin
(Of course the ``8" in \eqref{dandW} is just the value we happen to have below.)

\mn
{\em Proof.}
The proofs of (a) and (b) are essentially identical and we just give the former.
Set $W=\{z\in V\sm \{x\}:d_\A(z)\geq 2\}$.
Since the events $\{z\in B\}$ ($z\in V\sm \{x\}$) are negatively associated
(see Proposition~\ref{NAProp}),
Proposition~\ref{NAProp''} and the second condition in \eqref{dandW} give
\beq{PrBWs}
\Pr(|B\cap W|= s) ~\leq ~\Cc{\ww}{s}(k/n)^s
 ~<~ (\ww k/n)^s ~<~ n^{-(2\eps -o(1))s}.
\enq
On the other hand we assert that, with $\Q =\{B\cap A_i\neq\0 ~\forall i\in [d]\}$, we have
\beq{PrQB}
\Pr(\Q||B\cap W|=s) < \qqq ^{d-8s}.
\enq
To see this, condition on the value, $Z$, of $B\cap W$ (with $|Z|=s$), and let
\[
I= \{i\in [d]: B\cap A_i\cap W=\0\}.
\]
Then $|I|\geq d-8s$
(by the first condition in \eqref{dandW})
and $B$ must meet the members of $\{A_i:i\in I\}$ in $V\sm W$, where they are pairwise
disjoint.  By Proposition~\ref{q0prop},
$\Pr(B\cap (A_i\sm W)\neq\0|B\cap W=Z)<\qqq $ for each $i$.
But, given $\R_Z:=\{B\cap W=Z\}$,
$B\sm Z$ is a uniformly chosen
$(k-s)$-subset of $V\sm W$, so by Propositions~\ref{NAProp} and \ref{NAProp'}
the events $Q_i=\{B\cap (A_i\sm W)\neq\0\}$
are conditionally NA given $\R_Z$ (with $\Q=\cap_{i\in I}\Q_i$);
thus Proposition~\ref{NAProp''} gives
\[
\Pr(\Q|\R_Z)<\qqq ^{|I|}\leq \qqq ^{d-8s},
\]
which implies \eqref{PrQB}.

Finally, combining \eqref{PrBWs} and \eqref{PrQB}, we have
\begin{eqnarray*}
\Pr(\Q) &=& \sum_{s\geq 0}\Pr(|B\cap W|=s)\Pr(\Q||B\cap W|=s)\\
&< & \sum_{s\geq 0}n^{-(2\eps -o(1))s}\qqq ^{d-8s}\\
&=& \qqq ^d\sum_{s\geq 0}(n^{-(2\eps -o(1))}\qqq ^{-8})^s\sim \qqq ^d
\end{eqnarray*}\qed

\begin{cor}\label{ABCor}
Suppose either $\A$ is as in (a) of Proposition~\ref{ABProp} and
$\B $ is
chosen uniformly from the $b$-subsets of $\K\nex$,
or
$\A$ is as in (b) of the proposition and
$\B $ is
chosen uniformly from the $b$-subsets of $\{A\in \K:y\in A,x\not\in A\}$.
Then
\[
\Pr(B\cap A_i\neq\0 ~\forall B\in \B, i\in [d]) < (1+o(1))^b\qqq ^{db}.
\]
\end{cor}

\nin{\em Terminology.}  Recall that $\A$, $\B$ (two families of sets) are
{\em cross-intersecting} if
$A\cap B\neq \0~\forall A\in \A, B\in \B$.

\mn
{\em Proof.}  Again we just discuss the first case.
We may take $\B=\{B_1\dots B_b\}$ with
$B_i$ uniform from $\K\nex\sm \{B_1\dots B_{i-1}\}$.
Then, with $\Q_i=\{B_i\cap A_j\neq\0~\forall j\in [d]\}$, we have
\[
\Pr(\cap \Q_i)  \leq \prod\Pr(\Q_i) < (1+o(1))^b\qqq ^{db},
\]
with the second inequality given by Proposition~\ref{ABProp}.
(The first is obvious:
since the $B_i$'s are drawn without replacement, the probability that all are
drawn from those members of $\K\nex$ that meet all $A_j$'s is less than it would be
if they were drawn independently.)
\qed

\mn
{\em Proof of Lemma~\ref{ClaimA}.}
Let $\Q(x,r)$ be the event that there is some $\cee$ as in the lemma, with
$|\cee\nex| $ ($=|\cee|-d_\cee(x))$ $=r$, and let $\Q(x)=\cup_{r\geq 1}\Q(x,r)$.
By Proposition~\ref{HProp} it is enough to show that
(for any $x$)
\beq{Qx1n}
\Pr(\Q(x)\wedge\R)= o(1/n).
\enq
(Recall $\R$ was defined in the paragraph containing \eqref{dhxy} and \eqref{Wx}.)
Let
\[
\R_x = \{m\leq m_0; ~d(x)\leq\bbb; ~d(x,z)\leq 8 ~\forall z\in V\sm \{x\}; ~|W_x|\leq \ww\}.
\]
Then $\R_x\supseteq \R$, so for \eqref{Qx1n} it will be enough to bound
\[
\Pr(\Q(x)\wedge\R_x)\leq \sum_{r\geq 1}\Pr(\Q(x,r)\wedge\R_x).
\]

Set
\[
\sss(x,r)=\left\{\begin{array}{ll}
\{d(x)\geq \tau\} &\mbox{if $r \geq 2/\eps$,}\\
\{d(x)\geq \aaa-r\} &\mbox{if $r< 2/\eps$,}
\end{array}\right.
\]
and notice that $\sss(x,r)\supseteq \Q(x,r)$.
(For $r\geq 2/\eps$ this is contained in the definition of $\Q(x,r)$
(which promises $d_\cee(x)\geq \tau$), and for smaller $r$
it is given by $d(x)\geq d_\cee(x) =|\cee|-r \geq \aaa-r$.)
Thus we have
\begin{eqnarray}
\Pr(\Q(x,r)\wedge\R_x) &=&\Pr(\Q(x,r)\wedge\sss(x,r)\wedge\R_x) \nonumber\\
&\leq &\Pr(\sss(x,r))
\Pr(\Q(x,r)|\R_x\wedge \sss(x,r)).\label{QxrRx}
\end{eqnarray}

\medskip
For all but quite small $r$, a bound on the second factor in \eqref{QxrRx} will suffice for our purposes.
To bound this factor, we
condition on values $\h_x=\{A_1\dots A_d\}$
and $|\h\nex|=t$
satisfying $\sss(x,r)\wedge \R_x$
(in particular $d\leq \bbb$ and $t\leq m_0$);
thus $\h\nex $ is a uniform $t$-subset, say $\{B_1\dots B_t\}$, of $\K\nex $.
If $\Q(x,r)$ holds under this conditioning, then there are $I\sub [d]$ of size at least
$\tau$
and $J \sub[t]$ of size $r$
such that the families $\{A_i:i\in I\}$ and $\{B_j:j\in J\}$ are cross-intersecting
(namely,
each of the $r$ members of $\cee\nex$ meets each of the $d_\cee(x)\geq  \tau$
members of $\cee_x$).

The probability that this happens for a fixed $I$ and $J$ as above (note the remaining
randomization is in the choice of $B_j$'s) is, by Corollary~\ref{ABCor}, less than
$
(1+o(1))^r\qqq ^{\tau r}
$,
and it follows that the probability of $\Q(x,r)$ under the present conditioning---so also under
conditioning on $\sss(x,r)\wedge \R_x$---is less than
\begin{eqnarray}
\mbox{$\Cc{d}{\leq r}\Cc{m_0}{r}(1+o(1))^r\qqq ^{\tau r}$}
&<& \left[(1+o(1))\bbb m_0 n^{-(1-\eps)}\right]^r\nonumber\\
&<& n^{-(\eps-o(1))r}.\label{sumir}
\end{eqnarray}
Here the first factor on the left bounds
the number of possibilities for the $d-d_\cee(x)\leq r$
members of
$[d]\sm I$; the first inequality uses $d\leq \bbb$;
and the second uses
$\bbb m_0< (1+o(1))(\vp^*)^2n/k < n^{1-2\eps +o(1)}$
(see \eqref{vp*'}).

Thus, as suggested above, the second factor on the r.h.s. of \eqref{QxrRx} is enough for us unless $r$
is very small; namely,
\beq{sumxsumr}
\sum_{r>2/\eps}\Pr(\Q(x,r)|\R_x\wedge \sss(x,r)) = o(1/n).
\enq

For smaller $r$ we must use the factor
$\Pr(\sss(x,r))$ from \eqref{QxrRx} (together with \eqref{sumir}).
Here \eqref{tvp} gives
$\Pr(d_v=t)/\Pr(d_v=t+1)< n^{o(1)}$ for $t\in [\aaa-r,\aaa]$, which, since
$r< O(1)$, implies
\[
\Pr(\sss(x,r))< n^{o(1)}\Pr(d_x\geq \aaa+1) < n^{-1+o(1)}.
\]
Finally, recalling \eqref{sumir}, we find that
(for $r\leq 2/\eps$) the r.h.s. of \eqref{QxrRx} is less than
$n^{-1+o(1)}n^{-(\eps-o(1))r}= n^{-(1+r\eps-o(1))}$, yielding
\[
\sum_{r=1}^{\lfloor 2/\eps\rfloor}\Pr(\Q(x,r)\wedge\R_x) <
\sum_{r\geq 1}n^{-(1+r\eps-o(1))} = o(1/n),
\]
and combining this with \eqref{sumxsumr} gives \eqref{Qx1n}.
\qed

\section{Proof of Lemma \ref{ClaimB}}\label{PL2}

%Let $m_0$, $\ww$ and $\qqq $ be as in Section~\ref{PL1}.
%
We prove Lemma~\ref{ClaimB} in the following equivalent form.
\begin{lemma}\label{ClaimB'}
A.s. there do not exist $x,y\in V$ and $\f\sub \h_x$, $\g \sub \h_y$
with $|\f|=|\g |=\gl$ and $\f$, $\g $ cross-intersecting.
\end{lemma}
\nin

\nin
{\em Proof.}
Let $\Q(x,y)$ be the event described in Lemma~\ref{ClaimB'} and $\Q=\cup \Q(x,y)$.
We want $\Pr(\Q)=o(1)$, for which
it is enough to show that (for any $x,y$)
\beq{PrQxy}
\Pr(\Q(x,y)\wedge\R) < o(n^{-2}).
\enq

\medskip
For the proof of \eqref{PrQxy} we condition on
values of:  $m$ satisfying \eqref{mm} (so we may think of $\h$ as $\{A_i:i\in [m]\}$);
\[
\mbox{$I_x:=\{i\in [\hhh]:x\in A_i\}$, ~
$I_y:=\{i\in [\hhh]:y\in A_i\}$}
\]
 with $|I_x|,|I_y|\leq \bbb$ and
$|I_x\cap I_y|\leq 8$ (see \eqref{dhxy}); and a value of $(A_i:i\in I_x)$ for which
$|\{z\in V\sm \{x\}: |\{i:v\in A_i\}|\geq 2 \}| < \ww$ (see \eqref{Wx}).
If $\Q(x,y)$ holds (under this conditioning), then there are $J_x\sub I_x\sm I_y$ and $J_y\sub I_y\sm I_x$,
each of size $\gl-8$, with the families $\{A_i:i\in J_x\}$ and $\{A_j:j\in J_y\}$ cross-intersecting.

The probability that this happens for a given $J_x$, $(A_i:i\in J_x)$ and $J_y$ is,
by Corollary~\ref{ABCor}, at most $[(1+o(1))\qqq ^{\gl-8}]^{\gl-8}=\qqq ^{(1-o(1))\gl^2}$,
whence
\begin{eqnarray}
\Pr(\Q(x,y)\wedge\R)  &<&\Cc{\bbb}{\gl}^2 \qqq ^{(1-o(1))\gl^2}\nonumber\\
&<& \exp[\gl(2\log (e\gb/\gl) -(1-o(1))\gl\log(1/\qqq ))]\nonumber\\
&<& \exp[(1-o(1))\gl^2\log(1/q)]
<o(n^{-3}),\label{whencegl}
\end{eqnarray}
where the third inequality uses $\bbb <(1+o(1))\vp^*$ and $\gl \geq \sqrt{\log n}/\log(1/q)$
to say
$\log(e\bbb/\gl) =O(\log\log n)$ and the last uses
$\gl \geq 2\sqrt{\log n/\log(1/q)}$.\qed

\section{Proof of Lemma \ref{ClaimC}}\label{PL3}

Of our main points, Lemma~\ref{ClaimC} is the only one requiring
the full power of the assumption \eqref{avoid.generic}, as well as
the one requiring the most work:
there are several ways to handle it, but we so far don't see anything very compact.

We again
%(as  in the proof of Proposition~\ref{HProp})
condition on a value of $m$ satisfying \eqref{mm}
(so $\h$ is chosen uniformly from the $m$-subsets of $\K$), and then,
rather than dealing directly with $\h$,
find it easier to work with sets chosen
{\em independently} from $\K$, which makes
essentially no difference since $m$ is so small compared to $|\K|$.
Precisely, if $m$ satisfies \eqref{mm},
$B_1\dots B_m$ is a uniform $m$-subset of $\K$, $A_1\dots A_m$ are chosen
uniformly and independently from $\K$, and we set $\D=\{\mbox{$A_1\dots A_m$ are distinct}\}$, then for
any event $\B$ we have
\begin{eqnarray*}
\Pr(A_1\dots A_m\models \B)&\geq& \Pr(\D)\Pr(A_1\dots A_m\models \B|\D) \\
&=&\Pr(\D)\Pr(B_1\dots B_m\models \B),
\end{eqnarray*}
whence
\begin{eqnarray*}
\Pr(B_1\dots B_m\models \B) &\leq&
\Pr(A_1\dots A_m\models \B)/\Pr(\D) \\
&\leq &[1-m^2/\Cc{n}{k}]^{-1}\Pr(A_1\dots A_m\models \B)\\
&=&(1+o(1))\Pr(A_1\dots A_m\models \B).
\end{eqnarray*}
\medskip
It is thus enough to prove the following statement.
\begin{lemma}\label{ClaimC'}
Suppose $A_1\dots A_{\ghi}$ are drawn uniformly and independently from $\K$, and
let
$\Q $ be the event that
$\{A_1\dots A_{\ghi}\}$ is a clique
with
at most one vertex of degree greater than $\gl$ and none of degree greater than $\tau$.
Then
\[
\Pr(\Q) = o\left(\Cc{m}{\ghi}^{-1}\right).
\]
\end{lemma}

Given $A=(A_1\dots A_{\ghi})\in \K^{\ghi}$ we define several related quantities.
Write
$d_i(v)$ for the degree of $v$ in
the multiset $\{A_1\dots A_i\}$
and set $d_v=d_{\ghi}(v)$.
(We no longer default to $d_v=d_\h(v)$, since $\h$ plays no further role in this section.)
Note that we regard $A$ as given and sometimes (not always) suppress it in our notation; for example
$d_i(v)$ could also be written (say) $d_{A,i}(v)$.

We will need to distinguish two possibilities, depending on whether there
is or is not an $x$ with
$d_A(x)> \ddd$.
We treat these in parallel, the analysis in the second case eventually being
more or less contained in that for the first.
To this end we let $V'=V\sm \{x\}$ if we have specified such a high-degree $x$
and  $V'=V$ otherwise.

Set
$W_i=\{v\in V':d_i(v)=2\}$,
$Z_i=\{v\in V':d_i(v)\geq 3\}$,
$U_i=W_i\cup Z_i$, $W=W_{\ghi}$,
$Z =Z_{\ghi}$ and $U=U_{\ghi}$ ($=W\cup Z$).
% ($ = \{v:d_v\geq 2\}$).
In addition---now, for reasons which will appear below (see \eqref{A1}-\eqref{A4}), retaining $A$ in the notation---set

\[
s_i(A) = |A_i\cap W_{i-1}|,  ~~ r_i(A) = |A_i\cap Z_{i-1}|   ~~~\mbox{for $i\in [\ghi]$}
\]
(with $W_0=Z_0=\0$),
$\gs(A)= (s_1(A)\dots s_{\ghi}(A))$,
$\rho(A)= (r_1(A)\dots r_{\ghi}(A))$,
$
s(A)=\sum s_i(A)
$
and
$
r(A)=\sum r_i(A).
$
Notice that
\beq{sArA}
s(A)=|Z|
~~\mbox{and} ~~r(A)=\sum_{v\in Z}(d_v-3).
\enq

\medskip
Finally, set
\beq{Psi}
\Psi =\sum_{v\in Z}\left[\Cc{d_v}{2}-1\right]
\enq
and notice that
\[
\Psi
~=~2|Z| +\sum_{v\in Z}\left[\Cc{d_v}{2}-3\right]
~= ~2|Z| +\tfrac{1}{2}\sum_{v\in Z}(d_v-3)(d_v+2).
\]
We will only use this
when $d_v\leq \ddd$ for all $v\in V'$, in which case, in view of
\eqref{sArA}, we have
\beq{Xrs}
\Psi\leq 2s(A) +(\gl +2)r(A)/2 .
\enq

\bigskip
From this point we take $A=(A_1\dots A_{\ghi})$ with the $A_i$'s as in
Lemma~\ref{ClaimC'} (so chosen uniformly and independently from $\K$);
thus the quantities defined above ($d_i(v)$ through $\Psi$) become random variables
determined by $A$.
(Recall $\ww =\max\{\vp^2k^2/n, \log^6 n\}$
and $q=(1+2k^2\ww/(\qq n^2))\qq$
were defined in \eqref{ww}
and \eqref{qqq}.)

\begin{prop}\label{WZprop}
With probability $1-o\left(\C{m}{\ghi}^{-1}\right)$,

\mn
{\rm (a)}  $|U| < \ww~$  and

\mn
{\rm (b)}  $|Z| < \ghi/\eps ~=:\zz$.
\end{prop}
\nin

\mn
{\em Proof.}
Notice first that
\beq{target}
\Cc{m}{\ghi} ~<~ \exp[\ghi \log(em/\ghi)]~<~ \exp [(1/2+o(1))\ghi\log n],
\enq
since
$m/\ghi \leq 3m/\vp \sim 3n/k< n^{1/2+o(1)}$.

Since each $d_v$ has the binomial distribution ${\rm Bin}(\ghi,k/n)$,
we have (for all $v ,\ell$) $\Pr(d_v\geq \ell) < (k\ghi/n)^\ell/\ell !$,
whence
$\E|U| < k^2\ghi^2/(2n)$
and $\E|Z| < (k\ghi)^3/(6n^2)< n^{-(2\eps-o(1))}\ghi$.

On the other hand,
by Propositions~\ref{NAProp} and \ref{NAProp'},
the events $\{d_v\geq \ell\}$ are negatively associated for any $\ell$;
so the probabilities in question may be bounded using Corollary~\ref{CChern}.  For (b), we have
\beq{PrZ}
\Pr(|Z|> \zz) < n^{- (2\eps-o(1))\ghi/\eps} < n^{-\ghi}= o\left(\Cc{m}{\ghi}^{-1}\right).
\enq

The calculations for (a) are more annoying.
Here we set $k=\sqrt{\gz n}$ and $\mu=k^2\ghi^2/(2n)=\gz\ghi^2/2$
(our upper bound on $\E|U|$).  The desired inequality is
\[  %\beq{PrUw}
\mbox{$\Pr(|U|\geq \ww) = o\left(\C{m}{\ghi}^{-1}\right).$}
\]  %\enq
We first observe that this is true provided
\beq{gc1}
\ghi > 3\log n/\gz,
\enq
since then (using \eqref{Chern1} with $\gl =\mu$) we have ({\em cf.} \eqref{target})
\[
\Pr(|U|\geq \ww) \leq \Pr(|U|\geq 2\mu) < \exp[-\tfrac{3\mu}{8}]
= \exp[-\tfrac{3\gz\ghi^2}{16}]
<\exp[-\tfrac{9}{16}\ghi\log n].
\]
In particular \eqref{gc1} holds if (e.g.) $\gz\geq 2$, since then (according to \eqref{thinge})
we have $\ghi > (1-o(1))\frac{\log n}{-\log (1-e^{-\gz})}>3\log n/\gz$.
So we may assume
\[
\mbox{$\ghi\leq 3\log n/\gz~$ and $~\gz \leq 2$.}
\]
We then have $\log^6 n> 2\mu$, since
$\log^6n\leq 2\mu = \gz\ghi^2 \leq 9\log^2n/\gz$ implies $\gz <o(1)$,
yielding $\log(1/q) =\go(1)$ and $2\mu =\gz\ghi^2 = o((\vp^*)^2) =o(\log^6n)$, a contradiction.
Thus, again using Corollary~\ref{CChern}, we have
\[
\mbox{$\Pr(|U| >\ww)\leq \Pr(|U| >\log^6n) < \exp[-\gO(\log^6n)] < o\left(\C{m}{\ghi}^{-1}\right)$}
\]
(the last inequality holding since $\gz\leq 2$ implies $\ghi$ ($<\vp^*$) $=O(\log^3n)$).
\qed

\medskip
Set
%$\ww =k^2\ghi^2/n$, $\zz = \ghi/\eps$ and
\[
\sss = \{|W|\leq\ww , |Z|\leq \zz\}.
\]

\mn
By Proposition~\ref{WZprop},
Lemma~\ref{ClaimC'} will follow from

\beq{QRtoshow}
\Pr(\Q\wedge \sss) = o\left(\Cc{m}{\ghi}^{-1}\right).
\enq

\bigskip
For the proof of \eqref{QRtoshow} we will bound the probabilities
of various events whose union contains $\Q\wedge \sss$.
Set $\theta = \lfloor(n^\eps\log(1/q))^{-1}\rfloor$
and
\[
\A = \{\{A_1\dots A_{\ghi}\}~ \mbox{is a clique}\}.
\]
(Note $\theta$ need not be large---e.g. it will be zero for $k$ less than about
$\sqrt{\eps n\log n}$---so for once we do need the floor symbols.  The parts of
the following argument involving $\theta$ could be avoided when $\theta$ is small,
but there seems no point in treating this separately.)

For $x\in V$, $d\in (\ddd, \tau]$, and $\gs,\rho\in \Nn^{\ghi}$, let
\beq{A1}
\A(x,d,\rho,\gs)=  \A\wedge
\{d_x=d; ~d_v\leq \ddd ~\forall v\neq x; ~\rho(A)=\rho; ~\gs(A)=\gs\},
\enq
\beq{A2}
\A(x,d,\rho)=\A\wedge
\{d_x=d; ~d_v\leq \ddd ~\forall v\neq x; ~\rho(A)=\rho;~s(A)\leq \theta\},
\enq
\beq{A3}
\A(\rho,\gs)=\A\wedge
\{d_v\leq \ddd ~\forall v; ~\rho(A)=\rho;~\gs(A)=\gs\}
\enq
and
\beq{A4}
\A(\rho)=\A\wedge
\{d_v\leq \ddd ~\forall v;~\rho(A)=\rho; ~s(A)\leq \theta\}.
\enq

\mn
For $r,s\in \Nn$, let $X(r,s) = (\gl +2)r/2 +2s$
(the value in \eqref{Xrs}), and,
for $\varrho = (\varrho_1\dots \varrho_{\ghi})$, set $|\varrho|=\sum \varrho_i$.
\begin{lemma}\label{ML}
For any $x$, $d$, $\rho$, $\gs$ as above
with $|\rho|=r$ and $|\gs|=s$,
%$\rho =(r_1\dots r_{\ghi})$ and $\gs= (s_1\dots s_{\ghi})$ as above,
%with $\sum r_i=r$ and $\sum s_i=s$,
\beq{ML1}
\Pr(\A(x,d,\rho,\gs)\wedge\sss) < \C{\ghi}{d}\left(\frac{k}{n}\right)^d
\left(\frac{\zz k}{n}\right)^r\left(\frac{\ww k}{n}\right)^s \qqq ^{\Cc{\ghi}{2}-\Cc{d}{2} -X(r,s)}
\enq
and
\beq{ML2}
\Pr(\A(\rho,\gs)\wedge\sss) <
\left(\frac{\zz k}{n}\right)^r\left(\frac{\ww k}{n}\right)^s \qqq ^{\Cc{\ghi}{2} -X(r,s)}.
\enq
For any $x$, $d$, $\rho$ as above
with $|\rho|=r$,
\beq{ML3}
\Pr(\A(x,d,\rho)\wedge\sss) < \C{\ghi}{d}\left(\frac{k}{n}\right)^d
\left(\frac{\zz k}{n}\right)^r \qqq ^{\Cc{\ghi}{2}-\Cc{d}{2} -X(r,\theta)}
\enq
and
\beq{ML4}
\Pr(\A(\rho)\wedge\sss) <
\left(\frac{\zz k}{n}\right)^r \qqq ^{\Cc{\ghi}{2} -X(r,\theta)}.
\enq
\end{lemma}
\nin
(We will only use \eqref{ML1} and \eqref{ML2} with $s>\theta$.)

\medskip
Before proving Lemma~\ref{ML} we show that it implies \eqref{QRtoshow}.
Notice that $\Q$ is the (disjoint) union of the events
\beq{events}
\mbox{$\A(x,d,\rho,\gs),~
$
$\A(x,d,\rho), ~
$
$\A(\rho,\gs) ~
$
and
$~\A(\rho),
$}
\enq
where
$x\in  V$, $ d\in (\ddd,\tau]$, $\rho\in \mathbb N^{\ghi}$ and
$\gs\in \{(s_1\dots s_{\ghi})\in \Nn^{\ghi}: \sum s_i>\theta\}$.
Thus
\beq{QE}
\Pr(\Q\wedge\sss)\leq \sum \Pr(\eee\wedge\sss),
\enq
where $\eee$ ranges over the events in \eqref{events}.

It's now convenient to separate the contributions involving $x$, $\rho$
and $\gs$.  Set
\[
f(d) = n\C{\ghi}{d}\left(\frac{k}{n}\right)^d \qqq ^{-\C{d}{2}},
\]
\[
g(r) = \C{\ghi+r-1}{r}\left(\frac{\zz k}{n}\right)^r\qqq ^{-(\ddd+2)r/2},
\]
\[
h(s) = \C{\ghi+s-1}{s}\left(\frac{\ww k}{n}\right)^s\qqq ^{-2s}
\]
and
\[
h^* = \qqq ^{-2\theta}.
\]

Then, noting that (e.g.)
$|\{\rho\in \Nn^{\ghi}:|\rho| =r\}| = \Cc{\ghi+r-1}{r}$ and using
\eqref{ML1}-\eqref{ML4},
we find that $\Pr(\Q\wedge\sss)$ (or the r.h.s. of \eqref{QE})
is less than
\[
\qqq ^{\C{\ghi}{2}}\left[
\sum_{d,r,s} f(d)g(r)h(s) + h^*\sum_{d,r}  f(d)g(r) + \sum_{r,s}  g(r)h(s) +h^*\sum_r  g(r)\right],
\]
where
$ d$, $r$ and $s$ range over $(\ddd,\tau]$, $ \Nn$ and $(\theta,\infty)$ respectively.
Thus, since $\qqq ^\C{\ghi}{2}=o\left(\Cc{m}{\ghi}^{-1}\right)$ (by \eqref{thingb}, \eqref{m0m}
and \eqref{qqqqq} if $\ghi=\aaa$, and with plenty of room if
$\ghi = \vp^*/3$ ),
it is enough to show that each of
\[
\mbox{$\sum_{r\geq 0} g(r),$  $~~\sum_{s>\theta} h(s)~~$ and $ ~ ~h^*$}
\]
is $O(1)$ and that, with $F= \sum_{d\in (\ddd,\tau]}f(d)$,
\beq{gLF}
\qqq ^\C{\ghi}{2} F\Cc{m}{\ghi} ~~  (= \gL(\gc)F)  ~~ =o(1).
\enq
%(The first and third will be $o(1)$, the others $1+o(1)$.)
These are all easy calculations, as follows.

\medskip
First,
\[
g(r) \leq\left[ e\ghi (\zz k/n)n^{o(1)}\right]^r
< \left[ \ghi^2n^{-1/2+o(1)}\right]^r < n^{-(2\eps-o(1))r}
\]
(where the first inequality uses
$k > n^{1/2-o(1)} \Ra q> n^{-o(1)} \Ra \log(1/q) =o(\log n) \Ra \gl\log(1/q) = o(\log n)$),
implying $\sum_{r\geq 0} g(r) = 1+o(1)$.

\medskip
Second, since
\[
\Cc{\ghi+s-1}{s}^{1/s} < \tfrac{e(\ghi +s)}{s} < \tfrac{e(\ghi +\theta)}{\theta} < n^{\eps +o(1)}
\]
(for $s>\theta$), $\frac{\ww k}{n} < n^{-2\eps +o(1)} $ and $q=1-o(1)$, we have
\[
\mbox{$\sum_{s>\theta}h(s) < \sum_{s>\theta}n^{-(\eps -o(1))s} =o(1).$}
\]

\medskip
Third, $h^* =o(1)$ is immediate from our choice of $\theta$.

\medskip
The calculation for \eqref{gLF} requires a little more care.
Notice first that

\begin{eqnarray}
f(d)& <& n((e\ghi/d)(k/n))^d ~\qqq ^{-\C{d}{2}}\nonumber\\
&< & n\cdot n^{-(1/2 - o(1))d} \qqq ^{-\C{d}{2}}
 ~< ~ n\cdot \left[n^{-(1 - o(1))}\qqq ^{-d}\right]^{d/2} \label{f(d)}
\end{eqnarray}
(where the second inequality uses $\ghi/d <(1+o(1))\vp^*/\gl < n^{o(1)}$).
Here we may confine ourselves to
\beq{dbig}
d > (1-o(1))\log n/\log(1/q),
\enq
since for smaller $d$ the
expression in square brackets in \eqref{f(d)} is less than $n^{-\gO(1)}$
(and the exponent $d/2$ is at least $\gl/2=\go(1)$), so that the contribution of such $d$ to $F$
is $o(1)$.

For $d$ as in \eqref{dbig} the bound in \eqref{f(d)} is (rapidly) increasing in $d$ (passing from
$d$ to $d+1$ multiplies it by roughly $\sqrt{n}$);
so the contribution of such $d$ to $\gL(\gc)F$ is dominated by that of $d=\tau$.
For this term we have $\ghi = (1-\eps)^{-1}\tau = (1-\eps)^{-1}d > (1+\eps)\log n/\log(1/q)$ and

\begin{eqnarray*}
\gL(\ghi)f(\tau) &<& n^{-(1/2-o(1))\tau + \ghi/2 } q^{(\ghi-\tau)(\ghi +\tau-1)/2}\\
&<&n^{[\eps/2-(1+\eps)\eps (1-\eps/2) +o(1)]\ghi}\\
&=&
n^{-(\eps/2+\eps^2/2-\eps^3/2-o(1))\ghi}.
\end{eqnarray*}
Thus we have \eqref{gLF}.
\qed

\bigskip
For the proof of Lemma~\ref{ML}, we need the following easy observation.

\begin{prop}\label{martProp}
Let $Y_1\dots Y_\ell$ be
r.v.'s (not necessarily real-valued)
and write $y_i$ for a possible value of $Y_i$.
Let $\Z$ be a set of (``bad") prefixes $(y_1\dots y_i)$
closed under extension (i.e. $i<\ell$ and $(y_1\dots y_i)\in \Z$ imply
$(y_1\dots y_i,y_{i+1})\in \Z $ for every choice of $y_{i+1}$).
Set
\[
\Pr((Y_1\dots Y_i)\in \Z|y_1\dots y_{i-1}) = 1-\xi(y_1\dots y_{i-1}),
\]
where the conditioning has the obvious meaning and when $i=1$ the l.h.s.
is $\Pr((Y_1)\in \Z)$.
Then
\[
\Pr((Y_1\dots Y_\ell)\not\in Z)\leq
\max_{(y_1\dots y_{\ell})\not\in \Z}\prod_{i=1}^\ell \xi(y_1\dots y_{i-1})=:\xi.
\]
\end{prop}

\mn
{\em Proof.}
Define an auxiliary sequence $(X_0 \dots X_\ell)$ with $X_0\equiv 1$ and,
for $i\in [\ell]$,
\[
X_i= \left\{\begin{array}{ll}
0&\mbox{if $(Y_1\dots Y_i)\in \Z$,}\\
\xi(Y_1\dots Y_{i-1})^{-1}X_{i-1}&\mbox{otherwise.}
\end{array}\right.
\]

\mn
Then $\E X_\ell =X_0=1$
(since $(X_0 \dots X_\ell)$ is a martingale),
while
$X_\ell \geq \xi^{-1}$ whenever $(Y_1\dots Y_\ell)\not\in \Z$
(using the fact that $\Z$ is closed under extensions).  The conclusion follows.\qed

\medskip
We now turn to the proof of Lemma~\ref{ML},
beginning with the simpler \eqref{ML2} and \eqref{ML4};
the arguments for \eqref{ML1} and \eqref{ML3} are similar, and when we come to these
we will mainly just point out the necessary modifications.

For both \eqref{ML2} and \eqref{ML4} we will apply Proposition~\ref{martProp} to the
sequence $(Y_1\dots Y_{2\ghi})$,
where
\beq{Ys}
\mbox{$Y_{2j-1}=A_j\cap U_{j-1}~$ and $~Y_{2j} = A_j\sm U_{j-1}$.}
\enq
We first prove \eqref{ML2} and then discuss the changes needed for \eqref{ML4}.

\mn
{\em Proof of} \eqref{ML2}.
Here we say $(Y_1\dots Y_i)\in \Z$ (recall this is the set of ``bad" prefixes) if the associated $A_j$'s
(or parts of $A_j$'s)
satisfy at least one of:
\beq{Z1}
\mbox{$\{A_1\dots A_{\lfloor i/2\rfloor}\}$ is not a clique;}
\enq
\beq{Z2}
\mbox{for some $j\leq \lceil i/2\rceil$, $~|A_j\cap Z_{j-1}|\neq r_j~$ or $~|A_j\cap W_{j-1}|\neq s_j$;}
\enq
\beq{Z3}
\mbox{$|Z_{\lceil i/2\rceil}|> \zz$,  $~|W_{\lceil i/2\rceil}|> \ww~$ or
$d_{\lceil i/2\rceil}(v) > \ddd$ for some $v$}.
\enq

\mn
Then $\A(\rho,\gs)\wedge \sss=\{(Y_1\dots Y_{2\ghi})\not\in \Z\}$.

We next need to say something about the quantities
\[
\xi(y_1\dots y_{i-1}) = \Pr(Y_1\dots Y_i\not\in \Z|y_1\dots y_{i-1})
\]
appearing in Proposition~\ref{martProp},
where (we may assume)
$(y_1\dots y_{i-1})
\not\in \Z$.

If $i=2j-1$ then
\begin{eqnarray}\label{i=2j-1'}
\xi(y_1\dots y_{i-1}) &\leq & \Pr( |A_i\cap Z_{i-1}|\geq r_i, |A_i\cap W_{i-1}|\geq s_i|y_1\dots y_{i-1})\nonumber\\
&\leq & (\zz k/n)^{r_i}(\ww k/n)^{s_i}.
\end{eqnarray}
Here we again use
Propositions~\ref{NAProp} and \ref{NAProp''}, which,
since $(y_1\dots y_{i-1}) \not\in \Z$ implies $|Z_{i-1}|\leq \zz$ and $|W_{i-1}|\leq \ww$,
bound the probability in \eqref{i=2j-1'} by
\[
\Cc{\zz}{r_i} \Cc{\ww}{s_i}(k/n)^{r_i} ( k/n)^{s_i}.
\]

The case $i=2j$ is more interesting.  Here,
conditioning on the event
$
\{(Y_1\dots Y_{i-1})=(y_1\dots y_{i-1})\}
$,
we set
\beq{bj}
\bbb_j = \sum\{d_{j-1}(v):v\in A_j\cap U_{j-1}\}.
\enq
(Notice that this is determined by $(y_1\dots y_{i-1})$,
which includes specification of $Y_{2j-1}=A_j\cap U_{j-1}$.)
We will show
\beq{i=2j}
\xi(y_1\dots y_{i-1}) \leq  \qqq ^{j-1-\bbb_j}.
\enq

\medskip
Here we only consider \eqref{Z1}; that is, we
we ignore the requirements in \eqref{Z3}
(those in \eqref{Z2} are not affected by $Y_i$) and show that (given our conditioning)
the r.h.s. of \eqref{i=2j}
bounds the probability that $A_j$ meets all of $A_1\dots A_{j-1}$.
Now $A_j$ meets at most $\bbb_j$ members of $\{A_1\dots A_{j-1}\}$ in $U_{j-1}$, so to avoid
\eqref{Z1} must meet the $j-1-\bbb_j$ or more remaining members---say those
indexed by $I$---in $V\sm U_{j-1}$, where
they are pairwise disjoint.
This gives \eqref{i=2j}
since the events $Q_h=\{A_j\cap (A_h\sm U_{j-1})\neq\0\}$
($h\in I$) satisfy
$\Pr(Q_h) <\qqq $
(by Proposition~\ref{q0prop}) and are NA
(by Propositions~\ref{NAProp}
and \ref{NAProp'}), so by Proposition~\ref{NAProp''} we have
\[
\mbox{$\Pr(\cap_{h\in I} Q_h)\leq \prod_{h\in I} \Pr(Q_h) < \qqq ^{j-1-\bbb_j}$.}
\]

\bigskip
The last thing to notice here is that,
provided
$d_{\ghi}(v)\leq \ddd ~\forall v$---which in particular
is true whenever $(Y_1\dots Y_{2\ghi})\not\in \Z$; see \eqref{Z3})---we have
\beq{bjPsi}
\sum \bbb_j=\Psi\leq X(r,s)
\enq
(see \eqref{Psi} for $\Psi$ and \eqref{Xrs} for the inequality).
Finally, combining \eqref{i=2j-1'}, \eqref{i=2j} and \eqref{bjPsi}
(and $\sum_{j\in [\ghi]}(j-1)=\C{\ghi}{2}$) and applying
Proposition~\ref{martProp} gives \eqref{ML2}.\qed

\mn
{\em Proof of} \eqref{ML4}.
We now take $(Y_1\dots Y_i)\in \Z$ if the associated $A_j$'s
satisfy at least one of:
\beq{Z1A}
\mbox{$\{A_1\dots A_{\lfloor i/2\rfloor}\}$ is not a clique;}
\enq
\beq{Z2A}
\mbox{$\sum_{j\leq \lceil i/2\rceil}s_j(A)>\theta~$, or
for some $j\leq \lceil i/2\rceil$, $~|A_j\cap Z_{j-1}|\neq r_j$;}
\enq
\beq{Z3A}
\mbox{$|Z_{\lceil i/2\rceil}|> \zz$,  $~|W_{\lceil i/2\rceil}|> \ww~$ or
$d_{\lceil i/2\rceil}(v) > \ddd$ for some $v$}.
\enq

\mn
Then $\A(\rho)\wedge \sss\sub\{(Y_1\dots Y_{2\ghi})\not\in \Z\}$.

\medskip
The arguments bounding the quantities
\[
\xi(y_1\dots y_{i-1}) = \Pr(Y_1\dots Y_i\not\in \Z|y_1\dots y_{i-1})
\]
(again, for
$(y_1\dots y_{i-1})\not\in \Z$) are essentially identical to those above.
For $i=2j-1$ the bound
\beq{i=2j-1''}
\xi(y_1\dots y_{i-1}) \leq\Pr( |A_i\cap Z_{i-1}|\geq r_i|y_1\dots y_{i-1})
\leq (\zz k/n)^{r_i}
\enq
is justified in the same way as \eqref{i=2j-1'}.
For $i=2j$ we again define $\bbb_j$ as in \eqref{bj},
and \eqref{i=2j} follows as before.
(Note that our only reason for
retaining the constraint on $|W_{\lceil i/2\rceil}|$ in \eqref{Z3A}
is to enforce $\Pr(A_j\cap (A_h\sm U_{j-1})\neq\0) < \qqq $ in the proof of \eqref{i=2j}.

%that is, to rule out the (remote) possibility that $A_j$ and $A_h\sm U_{j-1}$
%are drawn from a universe of size significantly less than $n$.)

Finally,
\eqref{bjPsi} again holds
provided $(Y_1\dots Y_{2\ghi})\not\in \Z$
(this is where we use the first condition in \eqref{Z2A}), and
combining this with
\eqref{i=2j-1''} and \eqref{i=2j} we obtain
\eqref{ML4} {\em via}
Proposition~\ref{martProp}.\qed

\medskip
We now turn to the parts of Lemma~\ref{ML} involving $x$.
For $D\in \C{[\ghi]}{d}$ let
\[
\A(x,D,\rho,\gs)=  \A\wedge
\{x\in A_i\Leftrightarrow i\in D;
~d_v\leq \ddd ~\forall v\neq x; ~\rho(A)=\rho; ~\gs(A)=\gs\},
\]
\[
\A(x,D,\rho)=\A\wedge
\{x\in A_i\Leftrightarrow i\in D; ~d_v\leq \ddd ~\forall v\neq x; ~\rho(A)=\rho;~s(A)\leq \theta\}.
\]

\mn
Since
$
\Pr(\A(x,d,\rho,\gs) $ is the sum of the $ \Pr(\A(x,D,\rho,\gs))$'s
(and similarly for
$
\Pr(\A(x,d,\rho)$),
\eqref{ML1} and \eqref{ML3} will follow from (respectively)
\beq{ML1'}
\Pr(\A(x,D,\rho,\gs)) < \left(\frac{k}{n}\right)^d
\left(\frac{\zz k}{n}\right)^r\left(\frac{\ww k}{n}\right)^s \qqq ^{\Cc{t}{2}-\Cc{d}{2} -X(r,s)}
\enq
and
\beq{ML3'}
\Pr(\A(x,D,\rho)) < \left(\frac{k}{n}\right)^d
\left(\frac{\zz k}{n}\right)^r \qqq ^{\Cc{t}{2}-\Cc{d}{2} -X(r,\theta)}.
\enq

\nin
As the proofs of these closely track those of \eqref{ML2} and \eqref{ML4} (respectively),
with exactly the same modifications, we confine ourselves to
indicating what changes to the proof of \eqref{ML2} are needed for \eqref{ML1'}.

We again
apply Proposition~\ref{martProp}, in this case to the sequence $(Y_{1}\dots Y_{2\ghi})$ given by
\[
\mbox{$Y_{2j-1}=A_j\cap (U_{j-1}\cup\{x\})~$ and $~Y_{2j} = A_j\sm (U_{j-1}\cup\{x\})$}
\]
(which differs from \eqref{Ys} in the addition of $\{x\}$ to the $U_{j-1}$'s).
We say $(Y_{1}\dots Y_i)\in \Z$ if the associated $A_j$'s
satisfy at least one of
\eqref{Z1},
\eqref{Z2} (we recall for ease of reading that these were
\[
\mbox{$\{A_1\dots A_{\lfloor i/2\rfloor}\}$ is not a clique}
\]
and
\[
\mbox{for some $j\leq \lceil i/2\rceil$, $~|A_j\cap Z_{j-1}|\neq r_j~$ or $~|A_j\cap W_{j-1}|\neq s_j$}),
\]

\beq{Z3'}
\mbox{$|Z_{\lceil i/2\rceil}> \zz$,  $~|W_{\lceil i/2\rceil}|> \ww~$ or
$d_{\lceil i/2\rceil}(v) > \ddd$ for some $v\neq x$}
\enq
(which differs from \eqref{Z3} in the stipulation $v\neq x$) and
\beq{Z4B}
\mbox{for some $j\leq \lceil i/2\rceil$ either $j\in D$ and $x\not\in A_j$ or
$j\not\in D$ and $x\in A_j$.}
\enq

\mn
Then $\A(x,D,\rho,\gs)\wedge \sss=\{(Y_1\dots Y_{2\ghi})\not\in \Z\}$.

The bounds on the quantities
\[
\xi(y_1\dots y_{i-1})= \Pr(Y_1\dots Y_i\not\in \Z|y_1\dots y_{i-1})
\]
(again, for $(y_1\dots y_{i-1})  \not\in \Z$) are
modified as follows.
For $i=2j-1$ we use
\beq{i=2j-1B}
\xi(y_1\dots y_{i-1})
\leq \left\{\begin{array}{ll}
(k/n) (\zz k/n)^{r_i}(\ww k/n)^{s_i}&\mbox{if $j\in D$,}\\
(\zz k/n)^{r_i}(\ww k/n)^{s_i}&\mbox{otherwise;}
\end{array}\right.
\enq
this is justified ({\em via} Propositions~\ref{NAProp} and \ref{NAProp''})
in the same way as \eqref{i=2j-1'}.

\medskip
For $i=2j$
we define $\bbb_j$ as before
($\bbb_j = \sum\{d_{j-1}(v):v\in A_j\cap U_{j-1}\}$)
and set $\cc_j=[j-1]\sm D$
(again, a function of $(y_1\dots y_{i-1})$).
We then have
\beq{i=2jB}
\xi(y_1\dots y_{i-1}) \leq \left\{\begin{array}{ll}
\qqq ^{\cc_j-\bbb_j}&\mbox{if $j\in D$,}\\
\qqq ^{j-1-\bbb_j}&\mbox{otherwise.}
\end{array}\right.
\enq
The proof is essentially the same as that for
\eqref{i=2j}, the only difference being that when $j\in D$,
there is no requirement that $A_j$ meet those earlier $A_l$'s for which $l\in D$.
(On the other hand, the second bound in \eqref{i=2jB} uses the fact that
$x\not\in A_j$ (for $j\not\in D$), which follows from
$(y_1\dots y_{i-1})  \not\in \Z$; see \eqref{Z4B}.)

\medskip
Finally, applying Proposition~\ref{martProp} with
the combination of \eqref{i=2j-1B}, \eqref{i=2jB} and
$\sum \bbb_j=\Psi\leq X(r,s)$
(noted earlier in \eqref{bjPsi}) gives
\eqref{ML1'}, once we observe that
\[
\sum_{j\not\in D}(j-1) +\sum_{j\in D}\cc_j
=\sum_j(j-1) -\sum_{j\in D}|[j-1]\cap D|
= \C{\ghi}{2}-\C{d}{2}.
\]
\qed

\section{Large $\vp$}\label{Large}

\medskip
Here we
complete the proof of Theorem~\ref{MT} for $k$ as in \eqref{kbig}
by showing
\beq{largevp}
\mbox{for $\vp> \vp^*$, $\h$ satisfies EKR a.s.}
\enq
As already mentioned, this is mostly a matter of reducing to $\vp^*$ and applying
Lemmas~\ref{ClaimA}-\ref{ClaimC}.
(While there ought to be other ways to handle this, our main argument runs into
difficulties when $\vp$ is large, since the sets $W_x$, $W$, $Z$
used in the proofs of Lemmas~\ref{ClaimA}-\ref{ClaimC} are no longer small.)

\medskip
From now on we assume $\vp>\vp^*$.
We use the following natural reduction (coupling).
Setting $\rho = \vp^*/\vp$, we take $\h=\h_k(n,p)$ as usual and let $\g$ be
the random subhypergraph of $\h$ gotten by retaining edges independently, each with probability $\rho$;
thus $\g\sim \h_k(n,p^*)$, with $p^* = \vp^*/M$.

%We need two additional observations.
%
We would {\em like} to say that if EKR fails for $\h$, say at the nontrivial clique $\cee$,
then there is a decent chance
that the clique $\D:=\cee\cap \g$ fits one of the unlikely scenarios described in Lemmas~\ref{ClaimA}-\ref{ClaimC};
but this is not always true, since if $\cee$ is too close to trivial
then $\D$ is likely to actually {\em be} trivial.
This special situation is handled by the perhaps not uninteresting Lemma~\ref{L4},
and in other cases the desired reduction is given
by the routine Proposition~\ref{Reduction}.
%while the routine Proposition~\ref{Reduction} gives the desired reduction in less extreme cases.

\medskip
Set $r_0 = \xi \vp$ with $\xi = \log (1/q)/(2\log n)$
(as elsewhere, just a convenient value).
%\vp(5\gz e^{\gz})^{-1}$.

\begin{lemma}\label{L4}
A.s. there do not exist (in $\h$) a nontrivial clique $\cee$ and vertex $x$ such that
$|\cee|\geq \max\{\vp/2,d(x)\}$ and
$|\cee\nex|\leq r_0$.
\end{lemma}

\begin{prop}\label{Reduction}
Suppose $\cee$ is a nontrivial clique of $\h$ with $|\cee|\geq \gD\geq \vp/2$ and
$\gD_\cee\leq |\cee|-r_0$, and
let $x$ be a maximum degree vertex of $\cee$.
Then with probability at least $1/2-o(1)$, $\D:=\cee\cap \g$ satisfies:

\mn
{\rm (a)}  $|\D|\geq \max\{d_\g(x), \ghi\}$;

\mn
{\rm (b)}
$|\D\nex|> 2/\eps$;

\mn
{\rm (c)}
either $\gD_\D < \tau $ or $d_\D(x) > \gl$.
\end{prop}
\nin
(Recall $\ghi$, $\tau$ and $\gl$ were defined in \eqref{ghi}-\eqref{gl},
and note that in the present situation we have
$\gc=\vp^*/3$.)

\medskip
Before proving these assertions we show that they (with
Lemmas~\ref{ClaimA}-\ref{ClaimC})
give \eqref{largevp}.
Since $\gD\geq \vp/2$ a.s. (really, $d_v\geq \vp/2~ \forall v$ a.s. by Theorem~\ref{Chern}),
Lemma~\ref{L4} says it is enough to show that $\h$ is unlikely to contain a nontrivial clique $\cee$
with $|\cee|\geq \gD\geq \vp/2$ and $\gD_\cee < |\cee|-r_0$.
So we suppose this does happen, let $x$ be some maximum degree vertex of $\cee$,
and observe that
$\D$ and $x$ are then fairly likely (that is, with probability at least $1/2-o(1)$)
to exhibit one of the improbable behaviors described in
Lemmas~\ref{ClaimA}-\ref{ClaimC}; namely this is true if the conclusions of Proposition~\ref{Reduction} hold:

\mn
(i)  if $\D$ has at least two vertices of degree at least $\gl$, then Lemma~\ref{ClaimB} applies;
otherwise:

\mn
(ii) if $\gD_\D< \tau $ then we are in the situation of Lemma~\ref{ClaimC}
(since (a) of Proposition~\ref{Reduction}
gives $|\D|\geq \ghi$ and we assume
$\D$ has at most one vertex of degree at least $\gl$);

\mn
(iii)  if
$\gD_\D\geq  \tau $ then in fact $d_\D(x)\geq \tau $
(by (c), since we assume
$\D$ has at most one vertex of degree at least $\gl$ ($< \tau $));
so in view of (a) and (b)
we are in the situation described in Lemma~\ref{ClaimA}.
\qed

\mn
{\em Proof of Lemma~\ref{L4}.}
We need one preliminary observation.
For given $x$ and $\B\sub \K_x$, let $g(\B)$ be the probability that $A$ chosen uniformly from $\K\nex$
meets all members of $\B$.

Suppose that, for some $s$, $\B$ is a uniform $s$-subset of $\K_x$ and $A$ is
uniform from $\K\nex$ (these choices made independently).  Then
\beq{EgB}
\E g(\B) = \Pr(A\cap B\neq\0 ~\forall B\in\B) < \qq ^s,
\enq
the inequality holding because (i) $\Pr(A\cap B\neq \0) < \qq $ for $A$ and $B$ uniform from $\K\nex$ and $\K_x$ respectively,
and (ii) the probability in \eqref{EgB} is (obviously) no more than it would be if the members of $\B$
were chosen independently.
Markov's Inequality thus gives (for any $a\leq s$)
\[
\Pr(g(\B)> \qq ^a) < \qq ^{s-a}.
\]

\medskip
Now let
$
\sss = \pee\wedge\{d(x)\geq \vp/2 ~\forall x\}
$
(recall $\pee$ was defined in the paragraph containing \eqref{pee}),
noting that
(by \eqref{pee} and Theorem~\ref{Chern})
$\Pr(\ov{\sss})=o(1)$.
Let
\[
\Q(x)=
\{\exists \B\sub \h_x: ~ \mbox{$|\B|=d(x)-r_0$ and $g(\B)> \qq ^{\vp/4}$}\}
\]
and $\Q=\cup\Q(x)$.
Then
\begin{eqnarray}
\Pr(\Q\wedge\sss) &< &n\Cc{\bbb}{r_0}\qq ^{\vp/4-r_0}\nonumber\\
&<& n\exp[r_0\log (e\bbb/r_0)-(\vp/4-r_0)\log (1/\qq )].\label{QandR}
\end{eqnarray}
Recalling that $\vp\sim \bbb$
%and $\log (1/\qq )\sim e^{-\gz}$ (see \eqref{log1q}),
we have
\[
r_0\log (e\bbb/r_0)\sim \xi \vp\log (1/\xi) < (1/8)\vp \log (1/q)
\]
(since $\log(1/q) > n^{-1/4 +\gO(1)}$ implies $\log(1/\xi) < (1/4)\log n$);
so, noting that $q> n^{-o(1)}$ implies $r_0=o(\vp)$ and recalling
that $\vp> \vp^*$, we find that
the r.h.s. of \eqref{QandR} is $o(1)$.

Thus, with $\T$ the event in Lemma~\ref{L4}, the lemma will follow from
\beq{TQR'}
\Pr(\T\wedge\ov{\Q}\wedge\sss)=o(1).
\enq
We show
\beq{TQR}
\Pr(\T\wedge\ov{\Q}\wedge\sss)\leq n\sum_{r=1}^{ r_0}(\bbb m_0\qq ^{\vp/4})^r
\enq
(and then observe that the r.h.s. is small).

\mn
{\em Proof of} \eqref{TQR} {\em and} \eqref{TQR'} .
We consider occurrence of $\T$ at a given $x$, writing $\T(x)$ for this event.
Since
\[
\Pr(\T\wedge\ov{\Q}\wedge\sss)\leq \Pr(\T|\ov{\Q}\wedge\{d(x)\leq \bbb,m\leq m_0\})
\]
(the conditioning event contains $\ov{\Q}\wedge\sss$), it is enough to show
\[
\mbox{$\Pr^*(\T(x))< (\bbb m_0\qq ^{\vp/4})^r$,}
\]
where $\Pr^*$ denotes probability under conditioning on some
$\h_x$ of size at most $\bbb$ satisfying $\ov{\Q}(x)$,
together with a value $m\leq m_0$ of $|\h|$.

If, under this conditioning, $\T(x)$ occurs at $\cee$ with $|\cee\sm \h_x|=r$ ($\in [1,r_0]$),
then,
since $d_{\cee}(x)=|\cee|-r\geq \vp/2-r$ and $|\h_x\sm \cee|\leq r$,
there are $\B\sub \h_x$ and $\D\sub \h_{\ov{x}}$ (namely $\B=\cee_x$, $\D=\cee\nex$) with
\[
|\B|=|\cee|-r\geq \max\{d(x)-r,\vp/2-r\},
\]
$|\D|=r$, $\B$ and $\D$ cross-intersecting, and $g(\B)\leq \qq ^{\vp/4}$
(the last property implied by $\ov{\Q}(x)$;
of course if $|\B|\geq d(x)-r$ and $g(\B)> \qq ^{\vp/4}$,
then $g(\B')>\qq ^{\vp/4}$ for any $(d(x)-r_0$)-subset $\B'$ of $\B$).
But the probability that this occurs given $\h_x$ and $m$ as above
is at most
$\C{d(x)}{\leq r}\C{m-d(x)}{r}\qq ^{\vp r/4}< (\bbb m_0\qq ^{\vp/4})^r$
(which gives \eqref{TQR}).

\medskip
Finally (now for \eqref{TQR'}), we have
$
\bbb m_0\qq ^{\vp/4} < \vp^2n^{1/2 +o(1)}\qq ^{\vp/4} =o(1/n),
$
where the first inequality uses $\bbb\sim\vp$ and $m_0 \sim \vp n/k$,
and the second holds because $\vp^2\qq ^{\vp/4}$ is decreasing in $\vp>\vp^*$
and is $n^{-\go(1)}$ when $\vp=\vp^*$.\qed

\mn
{\em Proof of Proposition~\ref{Reduction}.}
Of course $\Pr(|\D|\geq d_\g(x))\geq 1/2$, so it's enough to show that each of the
other requirements (namely, $|\D|\geq \ghi$ and those in (b) and (c)) holds a.s.
These are all routine applications of Theorem~\ref{Chern} (or Corollary~\ref{CChern}):
First, $|\D|$ is binomial with mean
$|\cee|\rho \geq (\vp/2)\rho =\vp^*/2=3\ghi/2$, implying
$\Pr(|\D| <\ghi)< \exp[-\gO(\gc)]$.
Second, $\E |\D\nex|\geq r_0\rho=\xi\vp^* =\go(1)$,
so
$\Pr(|\D\nex|< 2/\eps ) < \exp[-\go(1)]$.
Third, since $\tau\gg\gl$ we have either
$\gD_\cee \rho $ ($= \E d_\D(x)$) $> 2\gl$, implying
$\Pr(d_\D(x) < \gl)=o(1)$, or
$\gD_\cee \rho <\tau/2$, implying
$\Pr(\gD_\D\geq \tau)< n\exp[-\gO(\tau)]=o(1)$;
thus (c) also holds a.s.\qed

\section{Small $k$}\label{Small}

Finally, we turn to the proof of Theorem~\ref{MT} for $k<n^{1/2-\gO(1)}$, say
\beq{klast}
k \leq n^{1/2-\eps}
\enq
with $\eps>0$ fixed
(note this is not the $\eps$ of Sections~\ref{Sketch}-\ref{Large}).
This is, as noted earlier, easier than what we've already done,
one reason being the absence of the issue
discussed following \eqref{chooset''}:
there will now always be an $\ga$ such that $\gD\geq \ga$ a.s.
and there is a.s. no nontrivial clique of size at least $\ga$.
This will mean that here we only need Proposition~\ref{DeltaProp}
(which for $k$ as in \eqref{klast} and {\em fixed} $\ga$ was proved in \cite{BBM})
and a simplified Lemma~\ref{ClaimC}.
Since most of this consists of easier versions of earlier arguments, parts of the
discussion will be a bit sketchy.

\medskip
It will be helpful to think of three regimes:
(i) $\vp <n^{-\gO(1)}$; (ii) $  n^{-o(1)}<\vp \ll 1$; and (iii) $\vp =\gO(1)$.
The last of these is treated in \cite[Theorem 1.1(iv)]{BBM}, so we may concentrate
on the first two.

We first need to specify $\ga$.  If we are in regime (ii) then we take $\ga$ as in
Section~\ref{Gen} (recall this assumed $\vp>n^{-o(1)}$ but not \eqref{kbig}),
noting that, in addition to $\gD\geq \ga$ a.s. (see \eqref{thingc}) and
$\gL(\ga) =o(1)$ (see \eqref{thingb}), we have $\ga =\go(1)$.
(Remark:  here $\ga=\ga_1$.)

In regime (i) we may assume (possibly passing to a subsequence of $n$'s) that
there is a (positive) integer $c$ such that
$n^{-1/c}\ll \vp =O(n^{-1/(c+1)}$;
we then take $\ga=c$, but will sometimes use $c$ to remind
ourselves that the value is a constant.
Here again we have
$\gD\geq \ga$ a.s.
(by Proposition~\ref{DeltaProp}
since $\Pr(d_v\geq \ga)\asymp \vp^\ga$),
as well as
$
\gL(\ga)=o(1),
$
which is given by \eqref{avoid.generic} once we observe that, by Harris' Inequality
\cite{Harris},
\begin{eqnarray}\label{Harris1}
\Pr(\gD\leq \ga) &=&\Pr(d_v\leq \ga ~\forall v)\nonumber\\
&\geq &\prod_v\Pr(d_v\leq \ga) = (1-O(1/n))^n = \gO(1).
\end{eqnarray}
(Of course if $\vp\ll n^{-1/(\ga +1)}$, then $\gD =\ga $ a.s., and it is not hard to see that if
$\vp \asymp n^{-1/(\ga+1)}$, then $\gD \in \{\ga ,\ga +1\}$ a.s. and each
possibility occurs with probability $\gO(1)$.)
Note also that we may assume $c ~(=\ga)~\geq 3$, since
if $c\leq 2$ then $\vp^2n\asymp \gL(2)=o(1)$ gives $\vp\ll n^{-1/2}$ and $\gD_\h \leq 1$ a.s.
(We may also note that if $c=3$ then $\gL(3)=o(1)$ implies $k\ll n^{1/3}$;
it is shown in
\cite{BBM} that for such $k$ EKR holds a.s. for {\em any} $\vp$.)

\medskip
In either regime we just need to show that $\h$ is unlikely to contain a nontrivial $\ga$-clique.
The arguments for the two regimes are similar and we treat them in parallel.
In each case we
will avoid some complications by first disposing of cliques with very large degrees
({\em cf.} Lemma~\ref{ClaimA}).

\medskip
If $\h$ contains a nontrivial clique of maximum degree at least $d$ then it contains a
``Hilton-Milner" family of size $d+1$; that is, $B_0\dots B_d$ such that $\cap_{i=1}^dB_i\sm B_0\neq\0$
and $B_i\cap B_0\neq 0~\forall i\in [d]$.
The probability that this occurs is, by Proposition~\ref{indeptAs}, within $o(1)$ of the
probability that it occurs for $A_1\dots A_m$ chosen independently from $\K$
(with $m$ chosen as usual).
The latter probability is less than
\beq{dlarge}
\Pr(m\not\hspace{-.01in}\models\eqref{mm})+\Cc{m_0}{d+1}(d+1)n(k/n)^d\qq^d
~<~ o(1)+ \vp ^{d+1}k^{2d-1}n^{-(d-2)}.
\enq
Here the factor $n$ on the l.h.s. is for a choice of $x\in \cap_{i=1}^dB_i\sm B_0$ and the inequality
uses $m_0\sim \mm =\vp n/k$ (with ``$\sim$" holding since, as already noted, we may assume $\vp=\gO(n^{-1/2})$
and, therefore, $\mm > n^{\gO(1)}$).
We then need to show that the r.h.s. of \eqref{dlarge} is $o(1)$ for suitable $d$.

For regime (i) we take $d =c-1$.
We have
\[
\gL(c)\asymp (\vp n/k)^c (k^2/n)^{\C{c}{2}} =
\left[\vp k^{c-2}n^{-(c-3)/2}\right]^{c},
\]
so $\gL(c)=o(1)$ implies
$
k^{c-2}\ll n^{(c-3)/2}/\vp.
$
Thus (for typographical reasons considering the $(c-2)^{{\rm nd}}$ power of the r.h.s. of \eqref{dlarge})
\begin{eqnarray*}
\left[\vp^c k^{2c-3}n^{-(c-3)}\right]^{c-2} &\ll&
\frac{\vp^{c(c-2)}n^{(2c-3)(c-3)/2}}{n^{(c-2)(c-3)}\vp^{2c-3}}
~=~ \left[\vp^{c-1}n^{1/2}\right]^{c-3}\\
&=& O(n^{-(c-1)/(c+1) +1/2})^{c-3} \\
&=& O(n^{-(c-3)^2/(2(c+1))}) ~=~ o(1),
\end{eqnarray*}
where we used $\vp = O(n^{-1/(c+1)})$ in the third step and $c\geq 3$ in the fourth.
Thus the r.h.s. of \eqref{dlarge} is $o(1)$.

For regime (ii) we take $d=\lfloor \ga/2\rfloor$ (say) and find that,
since $k< n^{1/2-\eps}$, the r.h.s. of \eqref{dlarge} is less than $n^{-\gO(\ga)}$.

\medskip
So (in either case) we just need to show that $\h$ is unlikely to contain a nontrivial $\ga$-clique
with maximum degree at most $d-1$ ($d$ as above).
The reduction to independent $A_i$'s preceding
Lemma ~\ref{ClaimC'} of course remains valid here, so the following analogue of
Lemma~\ref{ClaimC} completes the argument.

\begin{lemma}\label{ClaimC''}
Let $\ga$ be as above,
suppose $A_1\dots A_\ga$ are drawn uniformly and independently from $\K$, and
let
$\Q $ be the event that the multiset
$~\cee:=\{A_1\dots A_\ga\}$ is a nontrivial clique with $\gD_\cee \leq d-1$.
Then
$
\Pr(\Q) = o\left(\C{m_0}{\ga}^{-1}\right).
$
\end{lemma}

\mn
{\em Proof.}
This is a (much) simpler version of the proof of Lemma~\ref{ML}.
We retain the definitions of $d_i(v)$ and $d_v$ from that argument, but now set
$W_i=\{v:d_i(v)\geq 2\}$, $W=W_\ga$,
$s_i(A) = |A_i\cap W_{i-1}|,$ $\gs(A)=(s_1(A)\dots s_\ga(A))$,
\[
s(A)=\sum s_i(A)
=\sum_{v\in W}(d_v-2)
\]
and
\[   %\beq{Psi'}
\Psi =\sum_{v\in W}\left[\Cc{d_v}{2}-1\right] =
\tfrac{1}{2}\sum_{v\in W} (d_v+1)(d_v-2),
\]  %\enq
noting that if all $d_v$'s are at most $d$ then
\beq{Psi''}
\Psi\leq (d+1)s(A)/2.
\enq

For a counterpart of Proposition~\ref{WZprop}, with
$\ww = (\ga/\eps)$ ($\eps$ as in \eqref{klast}) and
$\sss = \{|W| \leq \ww\}$, we have
\[
\Pr(\ov{\sss})  = o(m_0^{-\ga})
\]
(since $\E|W|< (\vp k)^2/n < n^{-2\eps}$ implies
$\Pr(|W|\geq \ww) < n^{-2\eps \ww}= n^{-2\ga}\ll m_0^{-\ga}$).
So we need $\Pr(\Q \wedge\sss) =o\left(\C{m_0}{\ga}^{-1}\right)$.

\medskip
We again let $\A=\{\mbox{$\cee$ is a clique}\}$ and for
$\gs = (s_1\dots s_\ga)\in \mathbb N^\ga $ set
\[
\A(\gs) = \A\wedge \{\gD_\cee\leq d-1\}\wedge \sss \wedge \{\gs(A) = \gs\}.
\]
We have $\Q\wedge \sss =\cup_\gs \A(\gs)$ so, finally,
just need to show
\beq{Blast}
\mbox{$\sum_\gs \Pr(\A(\gs)) < o\left(\C{m_0}{\ga}^{-1}\right)$;}
\enq
with $q$ as in \eqref{qqq} (with the present $\ww$), this will follow from
\begin{lemma}\label{Llast}
For any $\gs$ as above with $\sum s_i=s$,
\beq{Bdlast}
\Pr(\A(\gs)) \leq \min\{(\ww k/n)^sq^{\C{\ga }{2}-(d+1)s/2}, (\ww k/n)^s\}.
\enq
\end{lemma}

\nin
Before sketching the proof of this, we show that it implies \eqref{Blast}, beginning with
regime (i) (so $\ga =c$, $d=c-2$ and $\C{m_0}{\ga}\asymp \mm^c$).
We use the first bound in \eqref{Bdlast} for $s:=|\gs|<c$ and the second for $s\geq c$.
For the latter we find that the contribution to $\mm^c\sum_{|\gs|\geq c}\Pr(\A(\gs))$
is at most
\[
\sum_{s\geq c}\Cc{s+c-1}{c-1}(\vp n/k)^c(\ww k/n)^s <
\sum_{s\geq c}((s+c)\vp\ww)^c (\ww k/n)^{s-c} = o(1).
\]

\nin
For the former,
the product of $\mm^c$ and the first bound in \eqref{Bdlast} is
\begin{eqnarray*}
(\vp n/k)^c(\ww k/n)^sq^{(c-1)(c-s)/2}
&\sim&
(\vp n/k)^c(\ww k/n)^s(k^2/n)^{(c-1)(c-s)/2}\\
&=& \vp^c \ww^s\left[(n/k)(k^2/n)^{(c-1)/2}\right]^{c-s}.
\end{eqnarray*}
If the expression in brackets is at most 1, then we have
\beq{mcsum}
\mm^c\sum_{|\gs|< c} \Pr(\A(\gs)) = O(\vp^c)
\enq
(since $\ww$ and $\C{s+c-1}{c-1}$ are $O(1)$, as is the number of terms in the sum),
and otherwise the sum in \eqref{mcsum} is on the order of
\[
\vp^c \left[(n/k)(k^2/n)^{(c-1)/2}\right]^c \asymp \gL(c) =o(1).
\]
%which we've assumed to be $o(1)$.

\medskip
For regime (ii), we use
the second bound in \eqref{Bdlast} for $s\geq 3\ga/2$, yielding
\begin{eqnarray*}
\Cc{m_0}{\ga}\sum_{|\gs|\geq 3\ga/2}\Pr(\A(\gs)) &<&
\sum_{s\geq 3\ga/2}\Cc{s+\ga -1}{\ga -1}(\vp n/k)^\ga (\ww k/n)^s \\
&<&
\sum_{s\geq 3\ga /2}((s+\ga )\vp\ww)^\ga  (\ww k/n)^{s-\ga } = o(1)
\end{eqnarray*}
(using $\ww=O(\ga)$ and $\ga\leq (1+o(1))\vp^* < n^{o(1)}$; see \eqref{thingf}, \eqref{vp*}).
On the other hand, the first bound in \eqref{Bdlast} gives
%(with $d$ as above)
\[    %\begin{eqnarray*}
\Cc{m_0}{\ga}\sum_{|\gs|<3\ga/2}\Pr(\A(\gs)) <
\Cc{m_0}{\ga}\sum_{s< 3\ga/2}\Cc{s+\ga -1}{\ga -1}(\ww k/n)^s q^{\C{\ga }{2}-(d+1)s/2}
= o(1)
\]   %\end{eqnarray*}
(because:  each of $\C{m_0}{\ga}$, $\Cc{s+\ga -1}{\ga -1}$
is
at most $\exp[O(\ga \log n)]$,
the $q$-term is less than $\exp[-\gO(\ga^2\log n)]$ (since $q<n^{-\gO(1)}$),
and, as noted above, $\ga=\go(1)$).

\qed

\mn
{\em Proof of Lemma~\ref{Llast}.}
This is
similar to the proof of Lemma~\ref{ML} and
we just indicate the little changes.
For the first bound in \eqref{Bdlast} we follow the proof of \eqref{ML2}
(beginning with the paragraph containing \eqref{Ys}),
with changes:
replace the $\ghi$'s by $c$'s and the $U$'s by $W$'s;
in \eqref{Z2} and \eqref{Z3} omit the condition involving $Z$
and replace $\gl$ by $d-1$ in \eqref{Z3};
omit the first factor in \eqref{i=2j-1'} (the proof doesn't change);
and replace $X(r,s)$ in \eqref{bjPsi} by $(d+1)s(A)/2$
(see \eqref{Psi''}).

For the second bound we use the same modifications and simply sacrifice
the contributions of the terms with $i=2j$ (so for these
we can just say $\xi(y_1\dots y_{i-1}) \leq 1$; thus the clique condition
\eqref{Z1} could be omitted here).\qed

\section{Necessity}\label{Nec}

Our main job in this section is to sketch the proof that the condition
\beq{star'}
\mbox{$\gL'(\gD) <o(1)$ a.s.}
\enq
in \eqref{NandS} is
necessary for EKR to hold a.s., but before doing so we say why Theorem~\ref{MT}
implies that it is sufficient.
Since $\h$ is a.s. EKR if either
\eqref{avoid.generic} holds or $\gD\leq 2$ a.s.
(the first by Theorem~\ref{MT}, the second trivially), we only need to consider
what happens when neither of these alternatives holds but \eqref{star'} does.
This means that each of $\{\gD\leq 2\}$ and $\{\gD\geq 3, ~\gL(\gD)<o(1)\}$ occurs
with probability $\gO(1)$ and their union occurs a.s., implying
$\vp\asymp n^{-1/3}$ and $\gL(3)<o(1)$.
But then $\gL(3)\asymp (\vp n/k)^3\qq^3$ implies $\qq=o(1)$,
so $k\ll \sqrt{n}$, $\qq\sim k^2/n$
and $\gL(3)\asymp k^3/n$.  Thus $k\ll n^{1/3}$ (again using $\gL(3)\ra 0$),
in which case EKR holds a.s. regardless of $\vp$
(as mentioned in Section~\ref{Small}, this was shown in \cite{BBM}; of course it
can also be extracted form the discussion in that section).

\medskip
We now turn to necessity.
We believe this actually holds for general $k$
(that is, without assuming \eqref{rangeofk}), but
our proof doesn't give this.  Of course, in view of the discussion preceding
Conjecture~\ref{EKRC}, the assertion seems less interesting for $k$
above about $\sqrt{(1/2)n\log n}$.

The proof of necessity
becomes easier (still not immediate)
if we retreat to, say, $k=O(\sqrt{n})$.
Here we give only a sketch of the argument, restricting to
$k$ as in \eqref{kbig} to avoid some annoyances,
with details---such as they are---mostly restricted to the more interesting points.
(Some instances of failure of EKR for smaller $k$ are given in \cite{BBM}.)

Note first of all that failure of \eqref{star'} means that there is some fixed $\gd>0$
such that, for infinitely many $n$,
$\Pr(\gL'(\gD)>\gd) >\gd$, whence also
\beq{gL'big}
\Pr(\gL(\gD)>\gd) >\gd;
\enq
so it is enough to show that \eqref{gL'big} implies
that EKR fails with probability at least some $\eta_\gd>0$.
(In what follows we just use $\gO(1)$'s.)

Set $\ga = \max \{t\in \mathbb N:\gL(t)>\gd\}$
and $\A=\{\gD\leq \ga\}$;
thus \eqref{gL'big} is
\beq{PrgD}
\Pr(\A) > \gd,
\enq
which we assume henceforth.

It is easy to see ({\em cf.} \eqref{thinge}) that
\beq{alphalast}
\ga \sim \frac{\log n}{\log (1/\qq)},
\enq
and we observe that (for any $v$)
\beq{a+1}
\Pr(d_v>\ga) =O(1/n),
\enq
since otherwise Proposition~\ref{DeltaProp} gives $\gD>\ga$ a.s.,
contradicting \eqref{PrgD}.

Here we do (finally) need some concrete notion of a ``generic" clique:
taking $\zz =\ga/\eps$
(with $\eps=1/4-c$ as in Sections~\ref{Sketch}-\ref{Large}),
%; {\em cf.} Proposition~\ref{WZprop}(b)),
say a clique---possibly with repeated edges---is
{\em generic} if it has maximum degree at most 3 and
at most $\zz$ vertices with degree equal to 3. Then with
\[
\B = \{\mbox{$\h $ contains a generic clique of size $\ga$}\},
\]
we will be done if we show %that (assuming \eqref{PrgD})
\beq{PrAB}
\Pr(\A\B) =\gO(1).
\enq
(The negative results of \cite{BBM} are achieved by showing (probable) existence of
$\gD$-cliques of maximum degree 2.)

Here again Proposition~\ref{indeptAs} allows us to work
with independent $A_i$'s; namely it implies that
\eqref{PrAB} will follow from:

\begin{lemma}\label{Lnec}
For any
$m$
satisfying \eqref{mm} and $\h=\{A_1\dots A_m\}$, with the $A_i$'s chosen
uniformly and independently from $\K$,
\beq{PPPAB}
\PPP(\A\B) =\gO(1).
\enq
\end{lemma}
\nin
(So we are using ``$\mathbb P$" for probabilities in this model.
Note $\h$ may now---in principle, though in reality essentially never---have repeated edges.)
%We now assume $m$, $\h$ and $A_1\dots A_m$ are as in Lemma~\ref{Lnec}
%and use ``$\mathbb P$" for probabilities based on these choices.

\medskip
We first assert that
\beq{PPPA}
\PPP(\A) =\gO(1).
\enq
This actually requires a little argument, but we just point out the difficulty.
That
\eqref{a+1}
implies the corresponding
$\PPP(d_v>\ga)=O(1/n)$ is easy,
the change in the distribution of $d_v$ from ${\rm Bin}(M,p)$ to ${\rm Bin}(m,k/n)$
having almost no effect.
But getting from this to \eqref{PPPA}---an implication which for
$\h_k(n,p)$ is given by Harris' Inequality; see \eqref{Harris1}---is no longer
immediate, since negative association now
works against us.

One way to handle this is to
compare the present $\h$ with $\h'=\h_k(n,p')$, with $p'>p$ chosen so that,
writing $\Pr'$ for the corresponding probabilities, we have
$\Pr'(d_v>\ga) =O(1/n)$ and $|\h'|\geq m$ a.s.
(We can then couple so that $\h'\supset\h$ a.s.---note $\h$ a.s. avoids
repeats---yielding $\PPP(\gD\leq \ga) > \Pr'(\gD\leq \ga) -o(1)=\gO(1)$.
Of course one must show there is such a $p'$, but
we omit this easy arithmetic.)

\medskip
For the proof of \eqref{PPPAB} we use the second moment method.
Set $N=[m]$ and
$\sss=\C{N}{\ga}$.
We now use
$\g$ for the set of generic $\ga$-cliques
(again, with repeated edges allowed).
For $S\sub N$ write
$A_S$ for the multiset $\{A_i:i\in S\}$ and $\gD_S$ for $\gD_{A_S}$
(so $\gD=\gD_N$).
In addition, set
$\B_S=\{A_S\in \g\}$, $X_S={\bf 1}_{\B_S}$ (these are
only of interest if  $S\in \sss$) and $X=\sum_{S\in \sss}X_S$.

We actually need estimates for the quantities $\E X_S$ and $\E X_SX_T$
(for $S,T\in \sss$) {\em conditioned on} $\A$, but will get these by
first dealing with the unconditional versions and then showing---the most
interesting point---that the
conditioning has little effect.  Thus we show (for any $S,T\in \sss$)
\beq{S1}
\E X_S \sim\qq^{\C{\ga}{2}};
\enq
\beq{S2}
\E X_SX_T < (1+o(1))\qq^{2\C{\ga}{2}-\C{|S\cap T|}{2}};
\enq
\beq{S3}
\mbox{$\E [X_S|\A]\sim\E X_S~~$ and $~~ \E [X_SX_T|\A]\sim\E X_SX_T$.}
\enq

We will say a little about the proofs of these main points below.
Once they are established we have, setting $\tilde{\E}[\cdot]=\E[\cdot|\A]$,
\[
\mu:=\tilde{\E} X \sim \Cc{m}{\ga}\qq^{\Cc{\ga}{2}} \sim \gL(\ga) =\gO(1)
\]
(using \eqref{m0m} for ``$\sim$") and an easy calculation gives
\begin{eqnarray*}
\tilde{\E} X^2 &=&\sum_S\sum_T \tilde{\E} X_SX_T\\
&<&
(1+o(1)) \Cc{m}{\ga}\qq^{2\Cc{\ga}{2}}\sum_{i=0}^\ga \Cc{\ga}{i}\Cc{m-\ga}{\ga-i}\qq^{-\Cc{i}{2}}
\sim \mu^2+\mu,
\end{eqnarray*}
whence
\[
\PPP(X\neq 0) \geq \mu^2/\tilde{\E} X^2 =\gO(1),
\]
which is what we want.

\medskip
The proofs of \eqref{S2}
and $\E X_S < (1+o(1))\qq^{\C{\ga}{2}}$ (for \eqref{S1})
are similar to (easier than) that
of Lemma~\ref{ClaimC'} and we will not pursue them here.
The proof of the reverse inequality in \eqref{S1} is also similar in spirit,
but less so in details.  We again think of choosing $A_1\dots A_\ga$ in order and
use $d_i$ for degrees in $\{A_1\dots A_i\}$.
Set $Z_i=\{v:d_i(v)\geq 3\}$,
$\Q_i=\{|Z_i|\leq \zz\}$,
$\R_i=\{A_i\cap Z_{i-1}=\0\}$,
$\T_i=\{A_i\cap A_j\neq\0~\forall j\in [i-1]\}$ and
\[
\B_i = \{\mbox{$\{A_1\dots A_i\}$ is a generic clique}\}.
\]
Then $\B_i=\B_{i-1}\R_i\T_i\Q_i$ and
\beq{PPPBi}
\PPP(\B_i) \geq \PPP(B_{i-1})\PPP(\R_i\T_i|B_{i-1}) -\PPP(\ov{\Q}_i).
\enq
We show by induction on $i$ (with $i=1$ trivial)
\beq{PPPBi'}
\PPP(\B_i) \geq (1-\gd_i)\qq^{\C{i}{2}},
\enq
for some $\gd_i < in^{-1/4+o(1)}$.
(This suffices because of \eqref{alphalast}, since $(\log(1/q))^{-1}< (1+o(1))n^c$;
see \eqref{log1q}.)

The relevant probabilities are bounded as follows.
First, the proof of Proposition~\ref{WZprop}(b) (see \eqref{PrZ}) gives
\beq{Qbd}
\PPP(\ov{\Q}_i) < \eta
\enq
for some $\eta < n^{-(2-o(1))\ga}$.
Second, trivially,
\beq{RBbd}
\PPP(\R_i|\B_{i-1}) \geq 1-\zz k/n
\enq
(this just uses $\B_{i-1}\sub \Q_{i-1}$).

\medskip
Third,
\beq{TRBbd}
\PPP(\T_i|\R_i\B_{i-1}) \geq (1-\theta)\qq^{i-1},
\enq
where $\theta < n^{-1/4+o(1)}$.
This one is less trivial than the first two.
We need the following general observation.

\begin{prop}\label{Pmeetall}
If $C_1\dots C_s, D_1\dots D_s$ are subsets of $V$
with $|C_i|=|D_i|$ $\forall i$ and the $D_i$'s pairwise
disjoint,
and $A$ is uniform from $\C{V}{k}$,
then
\beq{ACD}
\Pr(A\cap C_i\neq \0~\forall i)\geq\Pr(A\cap D_i\neq \0~\forall i).
\enq
\end{prop}

\nin
(This follows {\em via} induction from the fact---an easy coupling argument---that \eqref{ACD}
holds when $x\in C_i\cap C_j$ ($i\neq j$),
$D_i=C_i\sm\{x\}\cup \{y\}$ for some $y\in V\sm \cup C_\ell$, and $D_\ell = C_\ell$ for $\ell\neq i$.)

\medskip
By Proposition~\ref{Pmeetall}, the l.h.s. of \eqref{TRBbd} is
at least
\[
\Pr(A\cap A_j\neq\0 ~\forall j\in [i-1]),
\]
where
$A_1\dots A_{i-1}$ are (fixed) disjoint $(k-\zz$)-subsets of $U\in \C{V}{n-\zz}$ and $A $ is uniform
from $\C{U}{k}$.
Say $Y\sub U$ is {\em good} if $Y\cap A_j\neq\0 ~\forall j\in [i-1]$;
so we want
\beq{Agood}
\Pr(\mbox{$A$ is good}) \geq (1-\theta)\qq^{i-1}.
\enq
One way---there ought to be an easier one---to show this goes as follows.
Let $X\sub U$ be random with
each member of $U$ contained in $X$ with probability $\rho = (k-2\sqrt{k\ln n})/n$,
these choices made independently.
Then
\begin{eqnarray*}
\Pr(\mbox{$A$ is good}) &\geq &\Pr(\mbox{$X$ is good}) -\Pr(|X|>k)\\
&> &\Pr(\mbox{$X$ is good}) -n^{-2},
\end{eqnarray*}
where the first inequality holds because we can couple so that
$A\supseteq X$ whenever $|X|\leq k$, and the second is given by
Theorem~\ref{Chern}.
Thus, since $q^{i-1}> q^\ga > n^{-1-o(1)}$, \eqref{Agood} will
follow from
\beq{Xgood}
\Pr(\mbox{$X$ is good}) \geq (1-n^{-1/4+o(1)})\qq^{i-1}.
\enq
%This routine calculation is given in the appendix.\qed
%a calculation we omit.\qed
For verification of \eqref{Xgood},
set $\ell =k-\zz$.
Since
$\Pr(\mbox{$X$ is good}) = [1 - (1-\rho)^\ell]^{i-1}
$
and $i< \ga$, it is enough to show
\beq{ratio}
[(1 - (1-\rho)^\ell)/\qq]^\ga > 1-n^{-1/4+o(1)}.
\enq
Set (as earlier; see Section~\ref{Gen})
$\vt= 1-\qq = \tfrac{(n-k)_k}{(n)_k}
\sim e^{-k^2/n}
$
and define $\gc$ by
$
(1-\rho)^\ell = (1+\gc)\vt.
$
Then
$
(1 - (1-\rho)^\ell)/\qq =1- (\gc \vt/(1-\vt)),
$
so for \eqref{ratio} we need
$\ga \gc \vt/(1-\vt) < n^{-1/4+o(1)}$.
But it is easy to see
that we always have
$\ga\vt <O(\log n)$ (using \eqref{alphalast})
and $1-\vt $ ($=q$) $> n^{-o(1)}$
(since we assume \eqref{kbig});
so we really just need
\beq{gamma}
\gc < n^{-1/4+o(1)}.
\enq
Here we may expand
\begin{eqnarray*}
1+\gc=\frac{(1-\rho)^\ell}{\vt} &=&
\frac{(n-k)^k}{(n-k)_k}~\frac{(n)_k}{n^k}~(1-k/n)^{-\zz}
~\left(\frac{1-\rho}{1-k/n}\right)^{\ell}.
\end{eqnarray*}
The last two factors are at most
$1+k\zz/n +O(k^2\zz^2/n^2)< 1+n^{-1/4-\eps+o(1)}$ and $1+ O(n^{-1}k^{3/2}\sqrt{\log n}) < 1 + n^{-1/4+o(1)}$
(respectively),
while a little rearranging shows the product of the first two to be
\[
\prod_{j=0}^{k-1}\left( 1+ \frac{jk}{n(n-k-j)}\right) < 1+ k^3/n^2 < 1 + n^{-1/2+o(1)}.
\]
This proves \eqref{gamma} and finally establishes \eqref{TRBbd}.\qed

\medskip
By \eqref{Qbd}-\eqref{TRBbd} and \eqref{PPPBi'} for $i-1$, the r.h.s. of \eqref{PPPBi} is at least
\[
(1-\gd_{i-1})\qq^{\C{i-1}{2}}[(1-\zz k/n)(1-\theta) -\eta']\qq^{i-1}
~~~~~~~~~~~~~~~~~~~~
\]
\[
~~~~~~~~~~~~~~~~~~~~~~~~~~~~~>
(1-\gd_{i-1})[1-\{\zz k/n +\theta+\eta'\}]\qq^{\C{i}{2}},
\]
where we set $\eta'= \eta[(1-\gd_{i-1})\qq^{\C{i}{2}}]^{-1}$.
This gives \eqref{PPPBi'} since the expression in $\{~\}$'s
is
less than $n^{-1/4+o(1)}$.\qed

\medskip
Finally we turn to \eqref{S3}, for which we need the following observation.
\begin{prop}\label{Pasymp}
Let $s\in [m]$ and $t=m-s$.  Suppose
$S\in \C{N}{s}$
and $\D$ is an $s$-multisubset of $\K$ with $\gD_\D\leq C$.
If $\Pr(B(t,k/n) \geq \ga -C) =\rho/n$
then
\beq{PPPgD}
|\PPP(\gD\leq \ga|A_S=\D)-\PPP(\gD\leq \ga)|\leq sk\rho/n + n(sk/n)^{C+1}.
\enq
\end{prop}
\nin
{\em Proof.}
Let $B_S=\{B_i:i\in S\}$, where the $B_i$'s are chosen uniformly and independently
(of each other and the $A_j$'s) from $\K$,
and write $\gD^*$ for the maximum degree of $A_T\cup B_S$.
Set $V(X)=\{v:  d_X(v)>0\}$
(for $X$ a multisubset of $\K$).
On $\{A_S=\D\}$ we have
\beq{gD*gD}
\{\gD^*\leq \ga\}\sm \{\gD\leq \ga\}
~\sub ~ \{\exists v\in V(\D)~ d_T(v)> \ga -C\}
\enq
and
\beq{gDgD*}
\{\gD\leq \ga\}\sm \{\gD^*\leq \ga\}
~\sub ~ \{\gD_{B_S}> C\}\cup
\{\exists v\in V(B_S)~ d_T(v)> \ga -C\}
\enq
The probabilities of the event on the r.h.s. of \eqref{gD*gD}
and the second event on the r.h.s. of \eqref{gDgD*}
are at most $sk\rho/n$, and the probability of the first event on the r.h.s.
of \eqref{gDgD*} is less than $n(sk/n)^{C+1}$ (since $\E d_{B_S}(v)=sk/n)$.
The proposition follows.\qed

The arguments for the two statements in \eqref{S3} are nearly the same and we speak mainly
of the first.  This is equivalent to
$\PPP(\A|\B_S)\sim \PPP(\A)$ or, in view of
\eqref{PPPA},
$\PPP(\A|\B_S)=\PPP(\A)\pm o(1)$,
which will follow if we show that, for any generic $\ga$-clique $\D$,
\[
\PPP(\A|A_S=\D)=\PPP(\A)\pm o(1).
\]
This is, of course, an instance of Proposition~\ref{Pasymp}, for which we just have
to make sure that, with $s=\ga$ and $C=3$,
%and $\rho $ as in the proposition,
each part of the bound in
\eqref{PPPgD} is $o(1)$.
For the second part this is given by $\ga k/n < n^{-1/4-\eps +o(1)}$.
For the first, with $\xi=B(t,k/n)$, we have
$\Pr(\xi > \ga) \leq \Pr(B(m,k/n) >\ga) =O(1/n)$ (see \eqref{a+1}) and,
%({\em cf.} \eqref{tvp}),
for $u\sim \ga$,
\[
\frac{\Pr(\xi=u-1)}{\Pr(\xi=u)}=\frac{u(1-k/n)}{(t-u+1)k/n}\sim \frac{u}{tk/n}\sim\frac {\ga}{\vp}
< n^{o(1)}\]
(with the inequality given by \eqref{thingd}), whence $\rho <n^{o(1)}$
(and $sk\rho/n < n^{-1/4+o(1)}$).

(For the second part of \eqref{S3} we would have $s\in [\ga, 2\ga]$ and $C=6$.)\qed

\end{document}